\documentclass{amsart}
\usepackage{amsthm}
\newtheorem{thm}{Theorem}[section]
\newtheorem{cor}[thm]{Corollary}
\newtheorem{lem}[thm]{Lemma}

\newtheorem{prop}[thm]{Proposition}

\newtheorem{rem}[thm]{Remark}
\newtheorem{ex}[thm]{Example}

\newtheorem{defn}[thm]{Definition}

\numberwithin{equation}{section}
\newcommand{\R}{\mathbb{R}}
\newcommand{\N}{\mathbb{N}}
\begin{document}

\title[Viscosity solutions to second order PDEs on Riemannian manifolds]{
Viscosity solutions to second order partial differential equations
on Riemannian manifolds}
\author{Daniel Azagra, Juan Ferrera, Beatriz Sanz}

\address{Departamento de An{\'a}lisis Matem{\'a}tico\\ Facultad de
Matem{\'a}ticas\\ Universidad Complutense\\ 28040 Madrid, Spain}

\date{February 12, 2006}

\thanks{The authors were supported by grants MTM-2006-03531 and UCM-CAM-910626.}

\email{azagra@mat.ucm.es, ferrera@mat.ucm.es, bsanzalo@mat.ucm.es}



\begin{abstract}
We prove comparison, uniqueness and existence results for
viscosity solutions to a wide class of fully nonlinear second
order partial differential equations $F(x, u, du, d^{2}u)=0$
defined on a finite-dimensional Riemannian manifold $M$. Finest
results (with hypothesis that require the function $F$ to be
degenerate elliptic, that is nonincreasing in the second order
derivative variable, and uniformly continuous with respect to the
variable $x$) are obtained under the assumption that $M$ has
nonnegative sectional curvature, while, if one additionally
requires $F$ to depend on $d^{2}u$ in a uniformly continuous
manner, then comparison results are established with no
restrictive assumptions on curvature.
\end{abstract}

\maketitle

\section{Introduction}

The theory of viscosity solutions to nonlinear PDEs on
$\mathbb{R}^{n}$ (and on infinite-dimensional Banach spaces) was
introduced by M. G. Crandall and P. L. Lions in the 1980's. This
theory quickly gained popularity and was enriched and expanded
with numerous and important contributions from many
mathematicians. We cannot mention all of the significant papers in
the vast literature concerning viscosity solutions and
Hamilton-Jacobi equations, so we will content ourselves with
referring the reader to \cite{CLI} and the references given
therein.

More recently there have been various approaches to extend the
theory of viscosity solutions of first order Hamilton-Jacobi
equations, and the corresponding nonsmooth calculus, to the
setting of Riemannian manifolds. This is a natural thing to do,
because many functions arising from geometrical problems, such as
the distance function to a given set of a Riemannian manifold, are
not differentiable. Also, many important nonlinear equations full
of geometrical meaning, such as the eikonal equations, have no
classical solutions, and their {\em natural} solutions, which in
this case we think are the viscosity solutions, are not
differentiable (if some readers disagree with our saying that
viscosity solutions are the natural notion of solution for eikonal
equations, they might change their mind if they have a look at the
recent paper \cite{DevilleMatheron}, where the authors construct a
$1$-Lipschitz function $u$ defined on the closed unit ball
$\overline{B}$ of $\mathbb{R}^{n}, n\geq 2$, which is
differentiable on the open ball $B$, and such that $\|\nabla
u(x)\|=1$ almost everywhere, but $\nabla u(0)=0$; that is, the
eikonal equation $\|\nabla u(x)\|=1$ in $B$, $u=0$ on $\partial
B$, admits some exotic almost everywhere solutions which are
everywhere differentiable and are very different from its unique
viscosity solution, namely the distance function to the boundary
$\partial B$, which is not everywhere differentiable but is much
more natural from a geometric point of view).

Mantegazza and Menucci \cite{MantegazzaMenucci} studied viscosity
solutions to eikonal equations on Riemannian manifolds, in
connection with regularity properties of the distance function to
a compact subset of the manifold. In \cite{AFL2} a theory of
(first order) nonsmooth calculus for Riemannian manifolds
(possibly of infinite dimension) was introduced and applied to
show existence and uniqueness of viscosity solutions to
Hamilton-Jacobi equations on such manifolds. Simultaneously,
Ledyaev and Zhu \cite{LedyaevZhu} developed a (first order)
nonsmooth calculus on finite-dimensional Riemannian manifolds and
applied it to the study of Hamilton-Jacobi equations from a
somewhat different approach, related to control theory and
differential inclusions.

The usefulness of nonsmooth analysis on Riemannian manifolds has
been shown in \cite{GV}, where viscosity solutions are employed as
a technical tool to prove important results in conformal geometry.

However, to the best of our knowledge, no one has yet carried out
a systematic study of second order viscosity subdifferentials and
viscosity solutions to second order partial differential equations
on Riemannian manifolds.

In this paper we will initiate such a study by establishing
comparison, uniqueness and existence of viscosity solutions to
second order PDEs of the form
    $$
    F(x, u, du, d^{2}u)=0
    $$
where $u:M\to\mathbb{R}$ and $M$ is a finite-dimensional complete
Riemannian manifold. We will study the Dirichlet problem with a
simple boundary condition of the type $u=f$ on $\partial\Omega$,
where $\Omega$ is an open subset of $M$; and also the same
equation, with no boundary conditions, on all of $M$.

\medskip

Let us briefly describe the results of this paper. We begin with
the natural definition of second order subjet of a function
$u:M\to\mathbb{R}$, that is $J^{2,-}u(x)=\{(d\varphi(x),
d^{2}\varphi(x)) : \varphi\in C^{2}(M,\mathbb{R}), \,\,  f-\varphi
\textrm{ attains a local minimum at } x\}$. This is a nice
definition from a geometric point of view, but it would be
complicated and uneconomic to develop a nonsmooth calculus
exclusively based on this definition. It is more profitable to try
to localize the definition through charts and then use the second
order nonsmooth calculus on $\mathbb{R}^{n}$ to establish the
corresponding results on $M$. However, second derivatives of
composite functions are complicated, so not every chart serves
this purpose, and we have to work only with the exponential chart.
It is not difficult to see that $(\zeta, A)\in J^{2,-}u(x)$ if and
only if $(\zeta, A)\in J^{2,-}(u\circ\exp_{x})(0)$.

When one turns to the limiting subjet $\overline{J}^{2,-}u(x)$
(defined as the set of limits of sequences $(\zeta_{n}, A_{n})$,
where $(\zeta_{n}, A_{n})\in J^{2,-}u(x_{n})$ and $x_{n}$
converges to $x$), things become less obvious but, with the help
of a lemma which relates the second derivatives of a function
$\varphi:M\to\mathbb{R}$ to those of the function
$\psi=\varphi\circ\exp_{x}$ (at points near the origin in $TM_x$),
one can still show that $(\zeta, A)\in \overline{J}^{2,-}u(x)$ if
and only if $(\zeta, A)\in \overline{J}^{2,-}(u\circ\exp_{x})(0)$.

By using this characterization we can extend Theorem 3.2 of
\cite{CLI} to the Riemannian setting. This kind of result can be
regarded as a sophisticated nonsmooth fuzzy rule for the
superdifferential of the sum of two functions, and is the key to
the proof of all the comparison results in \cite{CLI} and in this
paper. The result essentially says that if $u_1$, $u_2$ are two
upper semicontinuous functions on $M$, $\varphi$ is a $C^2$ smooth
function on $M\times M$, and we assume that
$\omega(x_1,x_2)=u_1(x_1)+u_2(x_2)-\varphi(x_1, x_2)$ attains a
local maximum at $(\hat{x}_1,\hat{x}_2)$, then, for each
$\varepsilon>0$ there exist bilinear forms $B_i \in
\mathcal{L}^{2}_{s}((TM)_{\hat{x}_{i}}, \mathbb{R})$,  $i=1, 2$,
such that
$$\left(\frac{\partial}{\partial x_{i}}\varphi(\hat{x}_{1}, \hat{x}_{2}),B_i \right)\in \overline{J}^{ \, 2, +}u_i(\hat{x}_i)$$
for $i=1,...,k$, and the block diagonal matrix with entries $B_i$
satisfies
    $$-\left({1\over\varepsilon}+\|A\|\right)I\leq
    \left(
\begin{array}{cc}
  B_1   &    0    \\
   0    &  B_2   \\
\end{array}\right)\leq A+\varepsilon A^2, \eqno(*)
    $$
where $A=d^2\varphi(\hat{x}_{1}, \hat{x}_{2})\in
\mathcal{L}^{2}_{s}(TM_{\hat{x}_{1}}\times
TM_{\hat{x}_{2}},\mathbb{R})$. This is all done in Section $2$ of
the paper.

In the case $M=\mathbb{R}^{n}$ this result is usually applied with
$\varphi(x,y)=\frac{\alpha}{2}\|x-y\|^{2}$, whose second order
derivative is given by the matrix
    $$
    \alpha \left(
\begin{array}{cc}
  I   &    -I    \\
  -I  &     I  \\
\end{array}\right).
    $$
When applied to vectors of the form $(v,v)$ in
$\mathbb{R}^{n}\times \mathbb{R}^{n}$ this derivative vanishes,
which allows one to derive from $(*)$ that $B_{1}\leq B_{2}$ (as
quadratic forms). This in turn provides a very general form of
comparison result for viscosity solutions of the equation $F(x, u,
du, d^{2}u)=0$ in which the continuous function $F$ is assumed to
be degenerate elliptic (that is nonincreasing in the variable
$d^{2}u$), strongly increasing in the variable $u$, and uniformly
continuous with respect to $x$.

The natural approach in the Riemannian setting is then to consider
$\varphi(x,y)=\frac{\alpha}{2} d(x,y)^{2}$, where $d$ is the
Riemannian distance in $M$. Two problems immediately arise. First,
the function $\varphi$ is not differentiable in general if the
points $x, y$ are not suitably close to each other. This is
unimportant because, in the proof of the main comparison result,
we only need $\varphi$ to be $C^2$ smooth on a ball of small
radius around a point $x_{0}$ which is the limit of two different
sequences $x_{\alpha}$ and $y_{\alpha}$, and we have to evaluate
$d^{2}\varphi$ at the points $(x_{\alpha}, y_{\alpha})$.

The second problem, however, is substantial. The second derivative
of the function $\varphi$ is a quadratic form defined on
$TM_{x}\times TM_{y}$, and what we would like is that, when
applied to a vector of the form $(v, L_{xy}v)$, where $L_{xy}$ is
the parallel transport from $TM_{x}$ to $TM_{y}$ along the unique
minimizing geodesic connecting $x$ to $y$, this derivative is less
than or equal to zero. This way condition $(*)$ would imply that
$L_{\hat{x}_{2}\hat{x}_{1}}(B_2)\leq B_{1}$, where
$L_{\hat{x}_{2}\hat{x}_{1}}(B_2)$ is the parallel transport of the
quadratic form $B_2$ from $TM_{\hat{x}_{2}}$ to $TM_{\hat{x}_{1}}$
along the unique minimizing geodesic connecting $\hat{x}_2$ to
$\hat{x}_{1}$, defined by
$$\langle L_{\hat{x}_{2}\hat{x}_{1}}(B_{2})v, v\rangle
:=\langle B_{2}(L_{\hat{x}_{1}\hat{x}_{2}}v),
L_{\hat{x}_{1}\hat{x}_{2}}v\rangle. $$ And therefore we should be
able to conclude that, if $F$ is continuous, strongly increasing in
the variable $u$, and degenerate elliptic (that is $F(x,r,\zeta,
B)\leq F(x,r,\zeta, A)$ whenever $A\leq B$), then a natural
extension to $\mathcal{X}:=\{(x,r, \zeta, A) : x\in M,
r\in\mathbb{R}, \zeta\in TM_x, A\in\mathcal{L}^{2}_{s}(TM_x)\}$ of
the notion of uniform continuity of $F(x, r, \zeta, A)$ with respect
to the variable $x$ (namely, that
    $$
    | F(y,r, L_{xy}\zeta, L_{xy}P)-F(x,r, \zeta, P)|\to 0 \,\,
    \textrm{ uniformly as }  y\to x,
    $$
which we abbreviate by saying that $F$ is {\em intrinsically
uniformly continuous with respect to $x$}) would be enough to show
that comparison holds.

However, as we will show in Section $3$, one has that
    $$
    d^{2}\varphi(x,y)(v, L_{xy}v)^{2}\leq 0
    $$
for all $v\in TM_{x}$ if and only if $M$ has nonnegative sectional
curvature. Therefore, with this choice of $\varphi$, one can get
results as sharp as those in $\mathbb{R}^{n}$ only when one deals
with manifolds of nonnegative curvature. Nevertheless, if the
sectional curvature $K$ of $M$ is bounded below, say $K\geq-
K_{0}$, then one can show that
    $$
    d^{2}\varphi(x,y)(v, L_{xy}v)^{2}\leq 2 K_{0}d(x,y)^{2} \|v\|^{2}
    $$
for all $v\in TM_{x}$, and by using this estimation it is possible
to deduce that, if one additionally assumes that $F$ satisfies a
certain uniform continuity assumption with respect to the
variables $x$ {\em and} $D^2 u$ of the kind ``for every
$\varepsilon>0$ there exists $\delta>0$ such that
$d(x,y)\leq\delta$ and $P-L_{yx}Q\leq \delta I $ imply $F(y, r,
L_{xy}\zeta, Q)-F(x, r, \zeta, P)\leq\varepsilon$ for all
$\zeta\in TM_{x}^{*}$, $P\in\mathcal{L}^{2}_{s}(TM_{x}), \,
Q\in\mathcal{L}^{2}_{s}(TM_{y})$, $r\in\mathbb{R}$", then the
comparison principle holds for the equation $F=0$ (either with the
boundary condition $u=0$ on $\partial\Omega$, or with the
assumption that $M$ has no boundary and the functions $u, v$ for
which one seeks comparison are bounded). This is all shown in
Sections $4$ and $5$.

In Section $6$ we see that Perron's method works perfectly well in
the Riemannian setting. For instance one can show existence of
viscosity solutions to the equation $u+G(x, du, d^{2}u)=0$ on
compact manifolds under the same continuity assumptions on $G$ as
those that we require for comparison.

In particular, we get the following: if $M$ is a compact manifold
and $G$ is degenerate elliptic and uniformly continuous in the
above sense, then there exists a unique viscosity solution of
$u+G(x, du, d^{2}u)=0$ on $M$. If one additionally assumes that
$M$ has nonnegative sectional curvature then the above uniform
continuity assumption can be relaxed: it is enough to require that
$G$ is {\em intrinsically uniformly continuous with respect to
$x$}, meaning that ``for every $\varepsilon>0$ there exists
$\delta>0$ such that $d(x,y)\leq\delta$ implies $G(y,L_{xy}\zeta,
L_{xy}P)-G(x,\zeta, P)\leq\varepsilon$ for all $\zeta\in
TM_{x}^{*}$, $P\in \mathcal{L}^{2}_{s}(TM_{x})$".

We end the paper by discussing the applicability of the above
theory to some particular examples of equations.

\medskip

The notation we use is standard. $M=(M,g)$ will always be a
finite-dimensional Riemannian manifold. The letters $X, Y, Z, V,
W$ will stand for smooth vector fields on the Riemannian manifold
$M$, and $\nabla_{Y}X$ will always denote the covariant derivative
of $X$ along $Y$. The Riemannian curvature of $M$ will be denoted
by $R$. Geodesics in $M$ will be denoted by $\gamma$, $\sigma$,
and their velocity fields by $\gamma', \sigma'$. If $X$ is a
vector field along $\gamma$ we will often denote $X'(t)=
\frac{D}{dt}X(t)=\nabla_{\gamma'(t)}X(t)$. Recall that $X$ is said
to be parallel along $\gamma$ if $X'(t)=0$ for all $t$. The
Riemannian distance in $M$ will always be denoted by $d(x,y)$
(defined as the infimum of the lengths of all curves joining $x$
to $y$ in $M$).

We will often identify (via the Riemannian metric) the tangent
space of $M$ at a point $x$, denoted by $TM_{x}$, with the
cotangent space at $x$, denoted by $TM_{x}^{*}$. The space of
bilinear forms on $TM_{x}$ (respectively symmetric bilinear forms)
will be denoted by $\mathcal{L}^{2}(TM_{x})$ or
$\mathcal{L}^{2}(TM_{x}, \mathbb{R})$ (resp.
$\mathcal{L}^{2}_{s}(TM_{x})$ or $\mathcal{L}^{2}_{s}(TM_{x},
\mathbb{R})$). Elements of $\mathcal{L}^{2}(TM_{x})$ will be
denoted by the letters $A, B, P, Q$, and those of $TM_{x}^{*}$ by
$\zeta, \eta$, etc. Also, we will denote by $T_{2, s}(M)$ the
tensor bundle of symmetric bilinear forms, that is
$$
T_{2, s}(M)=\bigcup_{x\in M}\mathcal{L}^{2}_{s}(TM_x, \mathbb{R}),
$$
and $T_{2, s}(M)_x = \mathcal{L}^{2}_{s}(TM_x, \mathbb{R})$.

We will make extensive use of the exponential mapping $\exp_{x}$
and the parallel translation along a geodesic $\gamma$ throughout
the paper, and of Jacobi fields along $\gamma$ only in Section
$3$. Recall that for every $x\in M$ there exists a mapping
$\exp_{x}$, defined on a neighborhood of $0$ in the tangent space
$TM_x$, and taking values in $M$, which is a local diffeomorphism
and maps straight line segments passing through $0$ onto geodesic
segments in $M$ passing through $x$. The exponential mapping also
induces a local diffeomorphism on the cotangent space
$TM_{x}^{*}$, via the identification given by the metric, that
will be denoted by $\exp_{x}$ as well.

On the other hand, for a minimizing geodesic $\gamma:[0, \ell]\to
M$ connecting $x$ to $y$ in $M$, and for a vector $v\in TM_{x}$
there is a unique parallel vector field $P$ along $\gamma$ such
that $P(0)=v$, this is called the parallel translation of $v$
along $\gamma$. The mapping $TM_{x}\ni v\mapsto P(\ell)\in TM_{y}$
is a linear isometry from $TM_{x}$ onto $TM_{y}$ which we will
denote by $L_{xy}$. Its inverse is of course $L_{yx}$. This
isometry naturally induces an isometry (which we will still denote
by $L_{xy}$), $TM^{*}_x \ni\zeta\mapsto L_{xy}\zeta\in TM^{*}_y$,
defined by
    $$
    \langle L_{xy}\zeta, v\rangle_y :=\langle \zeta, L_{yx}v\rangle_x.
    $$
Similarly, $L_{xy}$ induces an isometry $\mathcal{L}^{2}(TM_{x},
\mathbb{R})\ni A\to L_{xy}(A)\in \mathcal{L}^{2}(TM_{y},
\mathbb{R})$ defined by
    $$
    \langle L_{xy}(A)v, v\rangle_y :=\langle A(L_{yx}v), L_{yx}v\rangle_x.
    $$
By $i_{M}(x)$ we will denote the injectivity radius of $M$ at $x$,
that is the supremum of the radius $r$ of all balls $B(0_{x}, r)$
in $TM_{x}$ for which $\exp_{x}$ is a diffeomorphism from
$B(0_{x}, r)$ onto $B(x,r)$. Similarly, $i(M)$ will denote the
global injectivity radius of $M$, that is $i(M)=\inf\{i_{M}(x)
: x\in M\}$. Recall that the function $x\mapsto i_{M}(x)$ is
continuous. In particular, if $M$ is compact, we always have
$i(M)>0$.

For Jacobi fields and any other unexplained terms of Riemannian
geometry used in Section $3$, we refer the reader to \cite{dC,
Sakai}.

%%%%%%%%%%%%%%%%%%%%

\medskip

%%%%%%%%%%%%%%%%%%%%%%%%%%%%%%%%%%%%%%%%%%%%%%%%%%%%%%%%%%%%%%%%%%%%%%%%%%%%%%%%

\section{Second order viscosity subdifferentials on Riemannian manifolds}

Recall that the Hessian $D^{2}\varphi$ of a $C^2$ smooth function
$\varphi$ on a Riemannian manifold $M$ is defined by
    $$
    D^{2}\varphi(X,Y)=\langle\nabla_{X}\nabla\varphi, Y\rangle,
    $$
where $\nabla\varphi$ is the gradient of $\varphi$ and $X$, $Y$
are vector fields on $M$ (see \cite{Sakai}, page 31). The Hessian
is a symmetric tensor field of type $(0,2)$ and, for a point $p\in
M$, the value $D^{2}\varphi(X,Y)(p)$ only depends of $f$ and the
vectors $X(p), Y(p)\in TM_{p}$. So we can define the second
derivative of $\varphi$ at $p$ as the symmetric bilinear form
$d^{2}\varphi(p):TM_{p}\times TM_{p}\to\mathbb{R}$
    $$
    (v,w)\mapsto d^{2}\varphi(p)(v,w):=D^{2}\varphi(X,Y)(p),
    $$
where $X, Y$ are any vector fields such that $X(p)=v, Y(p)=w$. A
useful way to compute $d^{2}\varphi(p)(v,v)$ is to take a geodesic
$\gamma$ with $\gamma'(0)=v$ and calculate
    $$
    \frac{d^{2}}{dt^{2}} \varphi(\gamma(t))|_{t=0},
    $$
which equals $d^{2}\varphi(p)(v,v)$. We will often write
$d^{2}\varphi(p)(v)^{2}$ instead of $d^{2}\varphi(p)(v,v)$.
%%%%%%%%%%%%%%
\begin{defn}
{\em Let $M$ be a finite-dimensional Riemannian manifold, and
$f:M\to (-\infty, +\infty]$ a lower semicontinuous function. We
define the second order subjet of $f$ at a point $x\in M$ by
    $$
    J^{2,-}f(x)=\{(d\varphi(x), d^{2}\varphi(x)) \, : \,
    \varphi\in C^{2}(M, \mathbb{R}), \,  f-\varphi \textrm{
    attains a local minimum at } x\}.
    $$
If $(\zeta, A)\in J^{2,-}f(x)$, we will say that $\zeta$ is a
first order subdifferential of $f$, and $A$ is a second order
subdifferential of $f$ at $x$.

Similarly, for an upper semicontinuous function $g:M\to[-\infty,
+\infty)$, we define the second order superjet of $f$ at $x$ by
    $$
    J^{2,+}f(x)=\{(d\varphi(x), d^{2}\varphi(x)) \, : \,
    \varphi\in C^{2}(M, \mathbb{R}), \,  f-\varphi \textrm{
    attains a local maximum at } x\}.
    $$
Observe that $J^{2,-}f(x)$ and $J^{2,+}f(x)$ are subsets of
$TM^{*}_{x}\times\mathcal{L}^{2}_{s}(TM_{x}, \mathbb{R})$, where
$\mathcal{L}^{2}_{s}(TM_{x}, \mathbb{R})\equiv
\mathcal{L}^{2}_{s}(TM_{x})$ stands for the symmetric bilinear
forms on $TM_{x}$. It is also clear that
$J^{2,-}f(x)=-J^{2,+}(-f)(x)$, and that we obtain the same
definitions if we replace the condition``$\varphi\in
C^{2}(M,\mathbb{R})$" with ``$\varphi$ is $C^2$ smooth on a
neighborhood of $x$".}
\end{defn}
%%%%%%%%%%%%%%

\medskip

By using the fact that a lower semicontinuous function $f$ is
bounded below on a neighborhood $B$ of any point $x$ with
$f(x)<\infty$, one can easily find a function $\varphi\in
C^{2}(M,\mathbb{R})$ such that $\inf_{y\in\partial
B}(f-\varphi)(y)>f(x)$, hence $f-\varphi$ attains a local minimum
at some point $z\in B$, and $(d\varphi(z),  d^{2}\varphi(z))\in
J^{2,-}f(z)$. This shows that the set
    $$
    \{z\in M : J^{2,-}f(z)\neq\emptyset\}
    $$
is dense in the set $\{x\in M: f(x)<\infty\}$. A similar statement
is true of upper semicontinuous functions. Therefore, when dealing
with semicontinuous functions, one has lots of points where these
subjets or superjets are nonempty, that is lots of points of
second order sub- or super-differentiability.

In the sequel $M$ will always denote an $n$-dimensional Riemannian
manifold. We next state and prove several results for subjets
which also hold, with obvious modifications, for superjets.
%%%%%%%%%%%%%%%%
\begin{prop}
Let $f:M\to (-\infty, +\infty]$ be a lower semicontinuous
function. Let $\zeta\in TM^{*}_{x}, A\in
\mathcal{L}^{2}_{s}(TM_{x}, \mathbb{R})$, $x\in M$. The following
statements are equivalent:
\begin{enumerate}
\item $(\zeta, A)\in J^{2,-}f(x)$.
\item $f(\exp_{x}(v))\geq f(x)+\langle\zeta, v\rangle_{x}+\frac{1}{2}\langle Av, v\rangle_{x}+ o(\|v\|^{2}).$
\end{enumerate}
\end{prop}
%%%%%%%%%%%%%%%%
\begin{proof} $(1)\implies(2)$: If $(\zeta, A)\in J^{2,-}f(x)$, by
definition there exists $\varphi\in C^{2}(M,\mathbb{R})$ such that
$f-\varphi$ attains a local minimum at $x$ and
$\zeta=d\varphi(x)$, $A=d^{2}\varphi(x)$. We may obviously assume
that $\varphi(x)=f(x)$, so we have
    $$
    f(y)-\varphi(y)\geq 0
    $$
on a neighborhood of $x$. Let us consider the function
$h(v)=\varphi(\exp_{x}(v))$ defined on a neighborhood of $0_x$ in
$TM_{x}$. We have that
    $$
    h(v)=h(0)+\langle dh(0), v\rangle_{x}+\frac{1}{2}\langle d^{2}h(0)v, v\rangle_{x}+ o(\|v\|^{2}).
    $$
By taking $y=\exp_{x}(v)$ and combining this with the above
inequality we get
    $$
    f(\exp_{x}(v))\geq f(x)+\langle dh(0), v\rangle_{x}+\frac{1}{2}\langle d^{2}h(0)v, v\rangle_{x}+
    o(\|v\|^{2}),
    $$
so we only need to show that $\zeta=dh(0)$ and $A=d^{2}h(0)$. To
see this, let us fix $v\in TM_x$ and consider the geodesic
$\gamma(t)=\exp_{x}(tv)$ and the function $t\mapsto
\varphi(\gamma(t))=h(tv)$. We have that
    $$
    \frac{d}{dt}h(tv)=\langle d\varphi(\gamma(t)), \gamma'(t)\rangle,
    $$
and
    $$
    \frac{d^{2}}{dt^{2}}h(tv) =
    \langle
    d^{2}\varphi(\gamma(t))\gamma'(t), \gamma'(t)\rangle.
    $$
In particular, for $t=0$, we get
    $$
    dh(0)(v)=\frac{d}{dt}h(tv) |_{t=0}=\langle d\varphi(x), v\rangle=\langle\zeta, v\rangle,
    $$
that is $dh(0)=\zeta$; and also
    $$
    \langle d^{2}h(0)v, v\rangle=\frac{d^{2}}{dt^{2}}h(tv) |_{t=0}=
    \langle
    d^{2}\varphi(x)v, v\rangle=\langle
    Av, v\rangle,
    $$
that is $A=d^{2}h(0)$.

\medskip

\noindent $(2)\implies(1)$: Define $F(v)=f(\exp_{x}(v))$ for $v$
in a neighborhood of $0_x\in TM_x$. We have that
    $$
    F(v)\geq F(0)+\langle \zeta, v\rangle_{x}+\frac{1}{2}\langle Av, v\rangle_{x}+
    o(\|v\|^{2}).
    $$
The result we want to prove is known to be true in the case when
$M=\mathbb{R}^{n}$, so there exists $\psi:TM_x\to\mathbb{R}$ such
that $F-\psi$ attains a minimum at $0$ and $d\psi(0)=\zeta$,
$d^{2}\psi(0)=A$. Since minima are preserved by composition with
diffemorphisms, the function $\varphi:=\psi\circ\exp_{x}^{-1}$,
defined on an open neighborhood of $x\in M$, has the property that
$f-\varphi=(F-\psi)\circ\exp_{x}^{-1}$ attains a local minimum at
$x=\exp_{x}^{-1}(0)$. Moreover, according to $(1)\implies(2)$
above, we have that
    $$
    d\varphi(x)=d\psi(0), \textrm{ and }
    d^{2}\varphi(x)=d^{2}\psi(0),
    $$
so we get $d\varphi(x)=\zeta$ and $d^{2}\varphi(x)=A$. Finally, by
using smooth partitions of unity we can extend $\varphi$ from an
open neighborhood of $x$ to all of $M$.
\end{proof}
%%%%%%%%%%%%%%%

%%%%%%%%%%%%%%%
\begin{cor}\label{def equivalente por composicion con exp}
Let $f:M\to (-\infty, +\infty]$ be a lower semicontinuous
function, and consider $\zeta\in TM^{*}_{x}, A\in
\mathcal{L}^{2}_{s}(TM_{x}, \mathbb{R})$, $x\in M$. Then
    $$
    (\zeta, A)\in J^{2,-}f(x) \iff (\zeta, A)\in
    J^{2,-}(f\circ\exp_{x})(0_x).
    $$
\end{cor}
%%%%%%%%%%%%%%%

Making use of the above characterization, one can easily extend
many known properties of the sets $J^{2,-}f(x)$ and $J^{2,+}f(x)$
from the Euclidean to the Riemannian setting. For instance, one
can immediately see that $J^{2,-}f(x)$ and $J^{2,+}f(x)$ are
convex subsets of $TM^{*}_{x}\times\mathcal{L}^{2}_{s}(TM_{x})$.
They are not necessarily closed, but if one fixes a $\zeta\in
TM^{*}_{x}$ then the set $\{A: (\zeta, A)\in J^{2,-}f(x)\}$ is
closed. A useful property that also extends from Euclidean to
Riemannian is the following: if $\psi$ is $C^2$ smooth on a
neighborhood of $x$ then
    $$
    J^{2,-}(f-\psi)(x)=\{\left(\zeta-d\psi(x),\, A-d^{2}\psi(x)\right) \, : \,
    (\zeta, A)\in J^{2,-}f(x)\}.
    $$
One can also see that $f$ is twice differentiable at a point $x\in
M$ (in the sense that for some (unique) $\zeta\in TM^{*}_{x},
A\in\mathcal{L}^{2}_{s}(TM_{x},\mathbb{R})$ we have that
$f(\exp_{x}(v))=f(x)+\langle \zeta, v\rangle+\frac{1}{2}\langle
Av,v\rangle+o(\|v\|^{2})$ as $v\to 0$) if and only if
$J^{2,-}f(x)\cap J^{2,+}f(x)$ is nonempty (in which case
$J^{2,-}f(x)\cap J^{2,+}f(x)=\{(\zeta, A)\}$).

\medskip

Next we have to define the closures of these set-valued mappings.
Let us first recall that a sequence $(A_n)$ with
$A_{n}\in\mathcal{L}^{2}_{s}(TM_{x_{n}})$ is said to converge to
$A\in\mathcal{L}^{2}_{s}(TM_{x})$ provided $x_{n}$ converges to
$x$ in $M$ and for every vector field $V$ defined on an open
neighborhood of $x$ we have that $\langle A_{n}V(x_{n}),
V(x_{n})\rangle$ converges to $\langle AV(x), V(x)\rangle$. Since
we have $\langle A V, W\rangle=\frac{1}{2}\left(\langle A(V+W),
V+W\rangle-\langle A V, V\rangle-\langle A W, W\rangle\right)$, it
is clear that this is equivalent to saying that $\langle
A_{n}V(x_{n}), W(x_{n})\rangle$ converges to $\langle AV(x),
W(x)\rangle$ for all vector fields $V, W$ on a neighborhood of $x$
in $M$.

Similarly, a sequence $(\zeta_n)$ with $\zeta_{n}\in
TM^{*}_{x_{n}}$ converges to $\zeta$ provided that $x_{n}\to x$
and $\langle \zeta_{n}, V(x_{n})\rangle \to \langle\zeta,
V(x)\rangle$ for every vector field $V$ defined on an open
neighborhood of $x$.
%%%%%%%%%%%%%%%%%%%
\begin{rem}\label{convergence senses coincide in Euclidean spaces}
{\em It is not difficult to see that, if $M=\mathbb{R}^{n}$, then
$A_n$ (respectively $\zeta_{n}$) converges to $A$ (resp. $\zeta$)
in the above sense if and only if $\|A_{n}-A\|\to 0$ (resp.
$\|\zeta_{n}-\zeta\|\to 0$) in
$\mathcal{L}^{2}_{s}(\mathbb{R}^{n}, \mathbb{R})$ (resp. in
$\mathbb{R}^{n}$).

It is also worth noting that $\|A_{n}-A\|\to 0$ (resp.
$\|\zeta_{n}-\zeta\|\to 0$) in
$\mathcal{L}^{2}_{s}(\mathbb{R}^{n}, \mathbb{R})$ (resp. in
$\mathbb{R}^{n}$) if and only if $\langle A_{n}v, v\rangle \to
\langle Av, v\rangle$ (resp. $\langle \zeta_{n}, v\rangle \to
\langle \zeta, v\rangle$) for every $v\in\mathbb{R}^{n}$, that is
pointwise convergence is equivalent to uniform convergence on
bounded sets, as far as linear or bilinear maps on
$\mathbb{R}^{n}$ are concerned.}
\end{rem}
%%%%%%%%%%%%%%%%%%%

%%%%%%%%%%%%%%%%%%%
\begin{defn}
{\em Let $f$ be a lower semicontinuous function defined on a
Riemannian manifold $M$, and $x\in M$. We define
    $$
    \overline{J}^{2,-}f(x)=\{(\zeta, A)\in TM^{*}_{x}\times\mathcal{L}_{s}(TM_x)
     \, : \, \exists x_{n}\in M, \exists (\zeta_{n}, A_{n})\in
     J^{2,-}f(x_{n})$$  $$ s.t.
         \, (x_{n}, f(x_{n}), \zeta_{n}, A_{n})\to (x,f(x),\zeta, A)\},
     $$
and for an upper semicontinuous function $g$ on $M$ we define
$\overline{J}^{2,+}g(x)$ in an obvious similar way.}
\end{defn}
%%%%%%%%%%%%%%%%%%%
\begin{rem}
{\em According to Remark \ref{convergence senses coincide in
Euclidean spaces}, we have that, in the case $M=\mathbb{R}^{n}$,
the sets $\overline{J}^{2,-}g(x)$ and $\overline{J}^{2,+}g(x)$
coincide with the subjets and superjets defined in \cite{CLI}.}
\end{rem}
%%%%%%%%%%%%%%%%%%%

\medskip

In order to establish the analogue of Corollary \ref{def
equivalente por composicion con exp} for the closure
$\overline{J}^{2,-}g(x)$, we will use the following fact.
%%%%%%%%%%%%%%%%%%%
\begin{lem}\label{relationship through exp}
Let $\varphi:M\to\mathbb{R}$ be a $C^2$ smooth function, and
define $\psi=\varphi\circ\exp_{x}$ on a neighborhood of a point
$0\in TM_{x}$. Let $\widetilde{V}$ be a vector field  defined on a
neighborhood of $0$ in $TM_{x}$, and consider the vector field
defined by $V(y)=d\exp_{x}(w_{y})(\widetilde{V}(w_{y}))$ on a
neighborhood of $x$ in $M$, where $w_{y}:=\exp_{x}^{-1}(y)$, and
let
    $$
    \sigma_{y}(t)=\exp_{x}(w_{y}+t\widetilde{V}(w_{y})).
    $$
Then we have that
    $$
    D^{2}\psi(\widetilde{V},\widetilde{V})(w_{y})=D^{2}\varphi(V,V)(y)+\langle\nabla\varphi(y), \sigma_{y}''(0)\rangle.
    $$
\end{lem}
%%%%%%%%%%%%%%%%%%%
\noindent Observe that $\sigma_{x}''(0)=0$ so, when $y=x$, we
obtain
    $$
    d^{2}\psi(0)(v,v)=d^{2}\varphi(x)(v,v)
    $$
for every $v\in TM_{x}$.
%%%%%%%%%%%%%%%%%%%
\begin{proof}
Fix $y$ near $x$. We have that
    $$
    \frac{d}{dt}\psi(w_{y}+t\widetilde{V}(w_{y}))=\frac{d}{dt}\varphi(\sigma_{y}(t))=
    \langle \nabla\varphi(\sigma_{y}(t)), \sigma_{y}'(t)\rangle,
    $$
and
    $$
    \frac{d^{2}}{dt^{2}}\psi(w_{y}+t\widetilde{V}(w_{y}))=\frac{d^{2}}{dt^{2}}\varphi(\sigma_{y}(t))=
    \langle\nabla_{\sigma_{y}'(t)}\nabla\varphi(\sigma_{y}(t)), \sigma_{y}'(t)\rangle+
    \langle\nabla\varphi(\sigma_{y}(t)), \sigma_{y}''(t)\rangle.
    $$
Note that $\sigma_{y}'(0)=V(y)$, hence by taking $t=0$ we get the
equality in the statement. Observe that when $y=x$ the curve
$\sigma_{x}$ is a geodesic, so $\sigma_{x}''(0)=0$.
\end{proof}
%%%%%%%%%%%%%%%%%%%

%%%%%%%%%%%%%%%%%%%
\begin{prop}\label{closures of subjets}
Let $f:M\to (-\infty, +\infty]$ be a lower semicontinuous
function, and consider $\zeta\in TM^{*}_{x}, A\in
\mathcal{L}^{2}_{s}(TM_{x}, \mathbb{R})$, $x\in M$. Then
    $$
    (\zeta, A)\in \overline{J}^{2,-}f(x) \iff (\zeta, A)\in
    \overline{J}^{2,-}(f\circ\exp_{x})(0_x).
    $$
\end{prop}
%%%%%%%%%%%%%%%%%%%
\begin{proof}
$(\Rightarrow):$ If $(\zeta, A)\in \overline{J}^{2,-}f(x)$ there
exist $x_{n}\to x$ and $(\zeta_{n}, A_{n})\in {J}^{2,-}f(x_n)$ so
that $\zeta_n \to \zeta, \, A_n \to A, \, f(x_n)\to f(x)$. Take
$\varphi_n \in C^{2}(M)$ such that $f-\varphi_{n}$ attains a
minimum at $x_n$ and $\zeta_{n}=d\varphi_{n}(x_{n}), \, A_n=
d^{2}\varphi_{n}(x_{n})$. Define $\psi_n =\varphi_ n\circ\exp_{x}$
on a neighborhood of $0$ in $TM_x$, and $v_n
=\exp_{x}^{-1}(x_{n})$. It is clear that $f\circ\exp_{x}-\psi_n$
attains a minimum at $v_n$. We then have that $(d\psi_n (v_n),
d^{2}\psi_n (v_n))\in J^{2,-}(f\circ\exp_x )(v_n)$, and since
$v_n\to 0$ and $f\circ\exp_x (v_n)\to f(x)$, we only have to show
that $d\psi_n (v_n)\to \zeta$ and $d^{2}\psi_n (v_n)\to A$.

Take a vector field $\widetilde{V}$ on $TM_x$, and define a
corresponding vector field $V$ on a neighborhood of $x$ in $M$ by
    $$
    V(y)=d\exp_{x}(w_{y})(\widetilde{V}(w_{y})),
    $$
where $w_{y}=\exp_{x}^{-1}(y)$. We have that
    $$
    \langle d\psi_n (v_n), \widetilde{V}(v_{n})\rangle=
    \langle d\varphi_n (x_n)\circ d\exp_{x}(v_n), \widetilde{V}(v_{n})\rangle=
    \langle d\varphi_n (x_n), V({x_{n}})\rangle,
    $$
so we get
    $$
    \langle d\psi_n (v_n), \widetilde{V}(v_n)\rangle=
    \langle \zeta_n, V(x_n)\rangle \to \langle\zeta, V(x)\rangle=
    \langle\zeta, \widetilde{V}(0)\rangle,
    $$
which shows $d\psi_n (v_n)\to\zeta$. On the other hand, according
to the preceding Lemma, we also have that
    $$
    d^{2}\psi_{n}(v_{n})(\widetilde{V}(v_{n}),\widetilde{V}(v_{n}))=
    A_{n}(V(x_{n}),V(x_{n}))+\langle \zeta_{n},
    \sigma_{x_{n}}''(0)\rangle,
    $$
where $\sigma_{y}(t)=\exp_{x}(w_{y}+t\widetilde{V}(w_{y}))$.

Notice that the mapping $y\mapsto \sigma_{y}''(0)$ defines a
smooth vector field on a neighborhood of $x$ in $M$ (and in
particular $\sigma_{x_{n}}''(0)\to \sigma_{x}''(0)=0$ as
$n\to\infty$). Since $A_{n}\to A$, $\zeta_{n}\to\zeta$, we get, by
taking limits as $n\to\infty$ in the above equality, that
$$d^{2}\psi_{n}(v_{n})(\widetilde{V}(v_{n}),\widetilde{V}(v_{n}))\to
A(V(x),V(x))+0 = A(\widetilde{V}(0),\widetilde{V}(0)),$$ which
proves that $d^{2}\psi_{n}(v_{n})\to A$.

\medskip

\noindent $(\Leftarrow)$ If $(\zeta, A)\in
\overline{J}^{2,-}(f\circ\exp_{x})(0)$ there exist $v_{n}\to 0$
and $(\widetilde{\zeta}_{n}, \widetilde{A}_{n})\in
{J}^{2,-}(f\circ\exp_{x})(v_n)$ so that $\widetilde{\zeta}_n \to
\zeta, \, \widetilde{A}_n \to A, \, f(x_n)\to f(x)$, where
$x_{n}=\exp_{x}(v_{n})$. Take $\psi_n \in C^{2}(TM_{x})$ such that
$f\circ\exp_{x}-\psi_{n}$ attains a minimum at $v_n$ and
$\widetilde{\zeta}_{n}=d\psi_{n}(v_{n}), \, \widetilde{A}_n=
d^{2}\psi(v_{n})$. Define $\varphi_n= \psi_{n}\circ\exp_{x}^{-1}$
on a neighborhood of $x$ in $M$. Then $f-\varphi_n$ attains a
minimum at $x_n$, so $(d\varphi_n (x_n), d^{2}\varphi_n (x_n))\in
J^{2,-}f(x_n)$, and we only have to show that $d\varphi_n (x_n)\to
\zeta$ and $d^{2}\varphi_n (x_n)\to A$. Take a vector field $V$ on
a neighborhood of $x$ in $M$, and define a corresponding vector
field $\widetilde{V}$ on a neighborhood of $0$ in $TM_x$ by
    $$
    \widetilde{V}(w_{y})=d\exp_{x}^{-1}(y)(V(y)),
    $$
where $w_{y}=\exp_{x}^{-1}(y)$. Now we have that
    $$
    \langle d\psi_n (v_n), \widetilde{V}(v_{n})\rangle=
    \langle d\varphi_n (x_n), V({x_{n}})\rangle,
    $$
from which we deduce that $d\varphi_n (x_n)\to\zeta$; and also, by
using this fact and the preceding Lemma,
\begin{eqnarray*}
& & d^{2}\varphi_{n}(V(x_{n}),V(x_{n}))=\\
& & \widetilde{A}_{n}(\widetilde{V}(v_{n}),\widetilde{V}(v_{n}))
    -\langle d\varphi_{n}(x_{n}), \sigma_{x_{n}}''(0)\rangle  \to\\
& & \to \langle AV(x), V(x)\rangle
    -\langle\zeta,\sigma_{x}''(0)\rangle=
    \langle AV(x), V(x)\rangle -0,
\end{eqnarray*}
concluding the proof.
\end{proof}
%%%%%%%%%%%%%%%%%%%

%%%%%%%%%%%%%%%%%%%
\begin{rem}
{\em One can see, as in the case of $J^{2, -}f(x)$, that if $\psi$
is $C^2$ smooth on a neighborhood of $x$ then
    $$
    \overline{J}^{\, 2,-}(f-\psi)(x)=\{\left(\zeta-d\psi(x),\, A-d^{2}\psi(x)\right) \, : \,
    (\zeta, A)\in \overline{J}^{\, 2,-}f(x)\}.
    $$
}
\end{rem}
%%%%%%%%%%%%%%%%%%%

\medskip

The following result is the Riemannian version of Theorem 3.2 in
\cite{CLI} and, as in that paper, will be the key to the proofs of
comparison and uniqueness results for viscosity solutions of
second order PDEs on Riemannian manifolds.
%%%%%%%%%%%%%%%%%%%
\begin{thm}\label{key to comparison}
Let $M_1 , ..., M_k$ be Riemannian manifolds, and $\Omega_i\subset
M_i$ open subsets. Define $\Omega=\Omega_1\times \ldots\times
\Omega_k\subset M_1\times\ldots\times M_k=M$. Let $u_i$ be upper
semicontinuous functions on $\Omega_i$, $i=1,...,k$; let $\varphi$
be a $C^2$ smooth function on $\Omega$ and set
    $$\omega(x)=u_1(x_1)+\ldots+u_k(x_k)$$
for $x=(x_1,...,x_k)\in\Omega$. Assume that
$\hat{x}=(\hat{x}_1,\ldots,\hat{x}_k)$ is a local maximum of
$\omega-\varphi$. Then, for each $\varepsilon>0$ there exist
bilinear forms $B_i \in
\mathcal{L}^{2}_{s}((TM_{i})_{\hat{x}_{i}}, \mathbb{R})$,  $i=1,
..., k$, such that
$$\left(\frac{\partial}{\partial x_{i}}\varphi(\hat{x}),B_i \right)\in \overline{J}^{\, 2, +}u_i(\hat{x}_i)$$
for $i=1,...,k$, and the block diagonal matrix with entries $B_i$
satisfies
    $$-\left({1\over\varepsilon}+\|A\|\right)I\leq
    \left(
\begin{array}{ccc}
  B_1   & \ldots &   0    \\
 \vdots & \ddots & \vdots \\
   0    & \ldots &  B_k   \\
\end{array}\right)\leq A+\varepsilon A^2,$$
where $A=d^2\varphi(\hat{x})\in
\mathcal{L}^{2}_{s}(TM_{\hat{x}},\mathbb{R})$.
\end{thm}
%%%%%%%%%%%%%%%%%%%
Recall that, for $\zeta\in TM^{*}$, $A\in\mathcal{L}(TM_{x}\times
TM_{x}, \mathbb{R})$, the norms $\|\zeta\|_{x}$ and $\|A\|_{x}$
are defined by
    $$
    \|\zeta\|_{x}=\sup\{ \langle\zeta, v\rangle_x \, : \, v\in TM_x, \|v\|_{x}\leq
    1\},
    $$
and
    $$
    \|A\|_{x}=\sup\{|\langle Av, v\rangle_x | \, : \, v\in TM_x, \|v\|_{x}\leq
    1\}=\sup\{ |\lambda| \, : \, \lambda \, \textrm{ is an
    eigenvalue of } A\}.
    $$
%%%%%%%%%%%%%%%%%%%%
\begin{proof}
The result is proved in \cite{CLI} in the case when all the
manifolds $M_i$ are Euclidean spaces, and we are going to reduce
the problem to this situation. By taking smaller neighborhoods of
the $x_i$ if necessary, we can assume that the $\Omega_i$ are
diffeomorphic images of balls by the exponential mappings
$\exp_{\hat{x}_{i}}:B(0, r_i)\to\Omega_i=B(\hat{x}_{i}, r_{i})$,
and that $\exp_{\hat{x}}$ maps diffeomorphically a ball in
$TM_{\hat{x}}$ onto a ball containing $\Omega$. The exponential
map $\exp_{\hat{x}}$ from this ball in
$TM_{\hat{x}}=(TM_{1})_{\hat{x}_1}\times \ldots \times
(TM_{k})_{\hat{x}_k}$ into $M$ is given by
    $$
    \exp_{\hat{x}}(v_{1}, ..., v_{k})=
    \left(\exp_{\hat{x}_1}(v_{1}), ..., \exp_{\hat{x}_k}(v_{k})\right).
    $$
Now define functions on open subsets of Euclidean spaces by
$\widetilde{\omega}(v)=\omega(\exp_{\hat{x}}(v))$ and
$\widetilde{u}_{i}(v_i)=u_{i}(\exp_{\hat{x}_i}(v_i))$. We have
that $\widetilde{\omega}(v_1,..., v_k)=\widetilde{u}_1
(v_1)+...+\widetilde{u}_k (v_k)$, and
$0_{\hat{x}}=(0_{\hat{x}_1},...,0_{\hat{x}_{k}})$ is a local
maximum of $\widetilde{\omega}-\psi$, where
$\psi=\varphi\circ\exp_{\hat{x}}$.

Then, by the known result for Euclidean spaces, for each
$\varepsilon>0$ there exist bilinear forms $B_i \in
\mathcal{L}^{2}_{s}((TM_{i})_{\hat{x}_{i}}, \mathbb{R})$,  $i=1,
..., k$, such that
$$\left(\frac{\partial}{\partial v_{i}}\psi(0_{\hat{x}}),B_i \right)\in \overline{J}^{\, 2, +}\widetilde{u}_i(0_{\hat{x}_i})$$
for $i=1,...,k$, and the block diagonal matrix with entries $B_i$
satisfies
    $$-\left({1\over\varepsilon}+\|A\|\right)I\leq
    \left(
\begin{array}{ccc}
  B_1   & \ldots &   0    \\
 \vdots & \ddots & \vdots \\
   0    & \ldots &  B_k   \\
\end{array}\right)\leq A+\varepsilon A^2,$$
where $A=d^2\psi(0_{\hat{x}})\in
\mathcal{L}^{2}_{s}(TM_{\hat{x}},\mathbb{R})$. According to
Proposition \ref{closures of subjets} we have that $$
\left(\frac{\partial}{\partial v_{i}}\psi(0_{\hat{x}}),B_i
\right)\in \overline{J}^{\, 2, +}\widetilde{u}_i(0_{\hat{x}_i}) \,
\, \iff \, \, \left(\frac{\partial}{\partial
v_{i}}\psi(0_{\hat{x}}),B_i \right)\in \overline{J}^{\, 2, +}
u_i(\hat{x}_i),
$$
so we are done if we only see that $$\frac{\partial}{\partial
v_{i}}\psi(0_{\hat{x}})=\frac{\partial}{\partial
x_{i}}\varphi(\hat{x}) \, \, \textrm{ and } \, \,
d^2\psi(0_{\hat{x}})=d^2\varphi(\hat{x}).$$ But this is a
consequence of Lemma \ref{relationship through exp}.
\end{proof}
%%%%%%%%%%%%%%%%%%%%

\medskip

Now we extend the notion of viscosity solution to a
Hamilton-Jacobi equation on a Riemannian manifold. In the sequel
we will denote
    $$
    \mathcal{X}:=\{(x,r, \zeta, A)
: x\in M, r\in\mathbb{R}, \zeta\in TM_x,
A\in\mathcal{L}^{2}_{s}(TM_x)\}
    $$

%%%%%%%%%%%%%%%%%%%
\begin{defn}[Viscosity solution]
{\em  Let $M$ be a Riemannian manifold, and
$F:\mathcal{X}\to\mathbb{R}$. We say that an upper semicontinuous
function $u:M\to\mathbb{R}$ is a viscosity subsolution of the
equation $F=0$ provided that
    $$
    F(x, u(x), \zeta, A)\leq 0
    $$
for all $x\in M$ and $(\zeta, A)\in J^{2,+}u(x)$. Similarly, a
viscosity supersolution of $F=0$ on $M$ is a lower semicontinuous
function $u:M\to\mathbb{R}$ such that
    $$
    F(x, u(x), \zeta, A)\geq 0
    $$
for every $x\in M$ and $(\zeta, A)\in J^{2,-}u(x)$. If $u$ is both
a viscosity subsolution and a viscosity supersolution of $F=0$, we
say that $u$ is a viscosity solution of $F=0$ on $M$.}
\end{defn}
%%%%%%%%%%%%%%%%%%%%%%%
\begin{rem}\label{continuity of F extends inequalities to closures of subjets}
{\em If $u$ is a solution of $F \leq 0$ and $F$ is continuous on
$\mathcal{X}$ then $F(x, u(x), \zeta, A)\leq 0$ for every $(\zeta,
A)\in \overline{J}^{\, 2,+}u(x)$. A similar observation applies to
solutions of $F\geq 0$ and solutions of $F=0$.}
\end{rem}
%%%%%%%%%%%%%%%%%%%%%%%

%%%%%%%%%%%%%%%%%%%%%%%
\begin{defn}[Degenerate ellipticity]\label{degenerate ellipticity}
{\em We will say that a function $F:\mathcal{X}\to\mathbb{R}$ is
{\em degenerate elliptic} provided that
    $$
    A\leq B \ \implies \, F(x,r,\zeta, B)\leq F(x,r,\zeta, A)
    $$
for all $x\in M, \, r\in\mathbb{R}, \, \zeta\in TM_{x}, \, A,
B\in\mathcal{L}^{2}_{s}(TM_{x})$.}
\end{defn}
%%%%%%%%%%%%%%%%%%%
\begin{ex}\label{examples of invariant functions}
{\em If we canonically identify the space of symmetric bilinear
forms on $TM_x$ with the space of self-adjoint linear mappings
from $TM_x$ into $TM_x$, we have that
    $$
    L_{yx}(Q)=L_{xy}^{-1}Q L_{xy}.
    $$
Hence
    $$
    \textrm{trace}(L_{yx}Q)=\textrm{trace}(Q), \,\, \textrm{ and } \,\,
    \textrm{det}_{+}(L_{yx}(Q))= \textrm{det}_{+}(Q)
    $$
(where $\textrm{det}_{+}A$ is defined as the product of the
nonnegative eigenvalues of $A$), and it is immediately seen that
the functions $G(r, \zeta, A)=-\textrm{det}_{+}(A)$ and
$H(r,\zeta, A)=-\textrm{trace}(A)$ are degenerate elliptic and,
moreover, are invariant by parallel translation, in the sense that
    $$
    G(r, \zeta, A)=G(r, L_{xy}\zeta, L_{xy}A)
    $$
for all $r\in\mathbb{R}, \zeta\in TM_x,
A\in\mathcal{L}^{2}_{s}(TM_x)$. The same can be said of all
nonincreasing, symmetric functions of the eigenvalues of $A$. Thus
one may combine such functions to construct many interesting
examples of equations to which our results apply, as we will see
later on.}
\end{ex}
%%%%%%%%%%%%%%%%%%%

%%%%%%%%%%%%%%%%%%%%%%%

%%%%%%%%%%%%%%%%%%%%%%%
\begin{rem}
If the function $F$ is degenerate elliptic, then every classical
solution of $F=0$ is a viscosity solution of $F=0$, {\em as is
immediately seen. However this may be not true if $F$ is not
degenerate elliptic; for instance when $M=\mathbb{R}$ the function
$u(x)=x^{2}-2$ is a classical solution of $u''(x)+u(x)-x^{2}=0$
but is not a viscosity solution.}
\end{rem}
%%%%%%%%%%%%%%%%%%%%%%%

In order that the theory of viscosity solutions applies to an
equation $F=0$, the following condition is usually required.

%%%%%%%%%%%%%%%%%%%%%%%
\begin{defn}[Properness]
{\em We will say that a function $F:\mathcal{X}\to\mathbb{R}$, $(x,
r, \zeta, A)\mapsto F(x, r, \zeta, A)$, is {\em proper} provided
\begin{enumerate}
\item[{(i)}] $F$ is degenerate elliptic, and
\item[{(ii)}] $F$ is nondecreasing in the variable $r$.
\end{enumerate}
}
\end{defn}
%%%%%%%%%%%%%%%%%%%%%%%

\medskip

%%%%%%%%%%%%%%%%%%%%%%%%%%%%%%%%%%%%%%%%%%%%%%%%%%%%%%%%%%%%%%%%%%%%%%%%%%%%%%%%%

\section{A key property of the Hessian of the function $(x,y)\mapsto
d(x,y)^{2}$}

When trying to establish comparison results for viscosity
solutions of second order PDEs on a Riemannian manifold $M$ we
will need to know how the Hessian of the function $\varphi:M\times
M\to\mathbb{R}$, $$\varphi(x,y)=d(x,y)^{2}$$ behaves. More
precisely we will need to know on which manifolds $M$ one has that
    $$
    d^{2}\varphi(x,y)(v, L_{xy}v)^{2}\leq 0
    \eqno(\sharp)
    $$
for all $v\in TM_{x}$, with $x,y\in M$ close enough to each other
so that $d(x,y)<\min\{ i_{M}(x), i_{M}(y)\}$.

Let us calculate this derivative. We have that
    $$
    \frac{\partial\varphi}{\partial x}(x,y)=2d(x,y)\frac{\partial d}{\partial
    x}(x,y)=-2\exp_{x}^{-1}(y). \eqno(1)
    $$
The second equality can be checked, for instance, by using the
first variation formula of the arc length (see \cite[p.
90]{Sakai}). Indeed, if $\alpha(t,s)$ is a variation through
geodesics of a minimizing geodesic $\gamma(t)$ with $y=\gamma(0)$
and $x=\gamma(\ell)$, where $\ell=d(x,y)$, and if $L(s)$ denotes
the length of the geodesic $t\mapsto \alpha(t,s)$, then
    $$
    \frac{d}{ds}L(s)|_{s=0}=\left[\langle V,T\rangle|_{0}^{\ell}-
    \int_{0}^{\ell}\langle V, \nabla_{T}T\rangle dt\right]=\left(\langle V(\ell), T(\ell)\rangle-
    \langle V(0), T(0)\rangle \right)
    $$
where $T=\partial\alpha/\partial t$ (so $\nabla_{T}T=0$) and
$V=\partial \alpha/\partial s$. Taking an $\alpha$ such that $V$
is the Jacobi field along $\gamma$ satisfying $V(0)=0$,
$V(\ell)=v$, we get
    $$
    \frac{\partial d}{\partial
x}(x,y)(v)=\frac{d}{ds}L(s)|_{s=0}=\frac{1}{\ell}\langle v,
-\exp_{x}^{-1}(y)\rangle.
    $$
Similarly, we have
    $$
    \frac{\partial\varphi}{\partial y}(x,y)=2d(x,y)\frac{\partial d}{\partial
    y}(x,y)=-2\exp_{y}^{-1}(x). \eqno(2)
    $$
Observe that $$\frac{\partial\varphi}{\partial
y}(x,y)+L_{xy}(\frac{\partial\varphi}{\partial x}(x,y)) = 0 =
\frac{\partial d}{\partial y}(x,y)+L_{xy}(\frac{\partial
d}{\partial x}(x,y)). \eqno(3)$$ By differentiating again in $(1)$
and $(2)$, we get
\begin{eqnarray*}
& &\frac{\partial^{2}\varphi}{\partial
x^{2}}(x,y)(v)^{2}=2\left(\frac{\partial d}{\partial
    x}(x,y)(v)\right)^{2}+2d(x,y)\frac{\partial^{2}d}{\partial
x^{2}}(x,y)(v)^{2}, \\
& &\frac{\partial^{2}\varphi}{\partial x\partial
y}(x,y)(v,w)=2\frac{\partial d}{\partial
    x}(x,y)(v)\frac{\partial d}{\partial
    y}(x,y)(w)+2d(x,y)\frac{\partial^{2}d}{\partial
x\partial y}(x,y)(v,w)\\
& &\frac{\partial^{2}\varphi}{\partial
y^{2}}(x,y)(w)^{2}=2\left(\frac{\partial d}{\partial
    y}(x,y)(w)\right)^{2}+2d(x,y)\frac{\partial^{2}d}{\partial
y^{2}}(x,y)(w)^{2},
\end{eqnarray*}
so, if we take $w=L_{xy}v$ and we sum the two first equations, and
then we use $(3)$, we get that
    $$
    \frac{\partial^{2}\varphi}{\partial
x^{2}}(x,y)(v)^{2}+\frac{\partial^{2}\varphi}{\partial x\partial
y}(x,y)(v,L_{xy}v)=2d(x,y)\left[ \frac{\partial^{2}d}{\partial
x^{2}}(x,y)(v)^{2}+\frac{\partial^{2}d}{\partial x\partial
y}(x,y)(v,L_{xy}v) \right]
    $$
and we get a similar equation by changing $x$ for $y$. By summing
these two equations we get
    $$
    d^{2}\varphi(x,y)(v, L_{xy}v)^{2}=
    2d(x,y) d^{2}(d)(x,y)(v, L_{xy}v)^{2},
    $$
so it is clear that condition $(\sharp)$ holds if and only if
    $$
    d^{2}(d)(x,y)(v, L_{xy}v)^{2}\leq 0
    \eqno(\flat)
    $$
for all $v\in TM_{x}$.

Another way to write conditions $(\sharp)$ or $(\flat)$ is
    $$
    \frac{d^{2}}{dt^{2}}\left(d(\sigma_{x}(t), \sigma_{y}(t))\right)|_{t=0} \leq
    0,\eqno(\natural)
    $$
where $\sigma_{x}$ and $\sigma_{y}$ are geodesics with
$\sigma_{x}(0)=x$, $\sigma_{y}(0)=y$, $\sigma_{x}'(0)=v$ and
$\sigma_{y}'(0)=L_{xy}v$. The function $t\mapsto h(t):=
d(\sigma_{x}(t), \sigma_{y}(t))$ measures the distance between the
geodesics $\sigma_{x}$ and $\sigma_{y}$ (which have the same
velocity and are parallel at $t=0$) evaluated at a point moving
along any of these geodesics.

We are going to show that the second derivative $h''(0)$ is
negative (that is, condition $(\natural)$ holds) if and only if
$M$ has positive sectional curvature.

In particular, by combining this fact with Equation $(3)$ (which
tells us that $h'(0)=0$), we see that the function $h(t)$ attains
a local maximum at $t=0$ if and only if $M$ has positive sectional
curvature. This corresponds to the intuitive notion that two
geodesics that are parallel at their starting points will get
closer if the sectional curvature is positive, while they will
spread apart if the sectional curvature is negative.
%%%%%%%%%%%%%%%%%%
\begin{prop}\label{key property of hessian of distance}
Condition $(\sharp)$ (equivalently $(\flat)$, or $(\natural)$)
holds for a Riemannian manifold $M$ if and only if $M$ has
nonnegative sectional curvature. In fact one has, for the function
$\varphi(x,y)=d(x,y)^{2}$ on $M\times M$, that:
\begin{enumerate}
\item If $M$ has nonnegative sectional curvature then
    $$
    d^{2}\varphi(x,y)(v, L_{xy}v)^{2}\leq 0
    $$
for all $v\in TM_{x}$, with $x,y\in M$ close enough to each other
so that $d(x,y)<\min\{ i_{M}(x), i_{M}(y)\}$.
\item If $M$ has nonpositive sectional curvature then
    $$
    d^{2}\varphi(x,y)(v, L_{xy}v)^{2}\geq 0
    $$
for all $v\in TM_{x}$, $x,y\in M$ such that $d(x,y)<\min\{
i_{M}(x), i_{M}(y)\}$.
\end{enumerate}
\end{prop}
%%%%%%%%%%%%%%%%%%
\noindent This fact must be known to the specialists in Riemannian
geometry, but we have been unable to find a reference for part
$(1)$, so we provide a proof. Let us begin by reviewing some
standard facts about the second variation of the arc length and
the energy functionals.

Take two points $x_{0}, y_{0}\in M$ with $d(x_{0},
y_{0})<\min\{i_{M}(x_{0}), i_{M}(y_{0})\}$, and let $\gamma$ be
the unique minimizing geodesic, parameterized by arc-length,
connecting $x_{0}$ to $y_{0}$. Denote $\ell=d(x_{0}, y_{0})$, the
length of $\gamma$. Consider $\alpha(t,s)$, a smooth variation of
$\gamma$, that is a smooth mapping $\alpha:[0,\ell]\times
[-\varepsilon, \varepsilon]\to M$ such that
$\alpha(t,0)=\gamma(t)$ for all $t\in [0,\ell]$. Consider the
length and the energy functionals, defined by
    $$
    L(s)=L(\alpha_{s})=\int_{0}^{\ell}\|\alpha_{s}'(t)\|dt
    $$
and
    $$
    E(s)=E(\alpha_{s})=\int_{0}^{\ell}\|\alpha_{s}'(t)\|^{2}dt,
    $$
where $\alpha_{s}$ is the variation curve defined by
$\alpha_{s}(t)=\alpha(t,s)$ for every $t\in [0,\ell]$. According
to the Cauchy-Schwarz inequality (applied to the functions
$f\equiv 1$ and $g(t)=\|\alpha_{s}'(t)\|$ on the interval $[0,
\ell ]$) we have that
    $$
    L(s)^{2}\leq\ell E(s),
    $$
with equality if and only if $\|\alpha_{s}'(t)\|$ is constant.
Therefore, in the case when $\alpha_{s}$ is a geodesic for each
$s$ (that is $\alpha$ is a variation of $\gamma$ through
geodesics) we have that
    $$
    L(s)^{2}=\ell E(s)
    $$
for every $s\in [-\varepsilon, \varepsilon]$.

Now take a vector $v\in TM_{x_{0}}$, set $w=L_{x_{0}y_{0}}v$, and
consider the geodesics $\sigma_{x_{0}}, \sigma_{y_{0}}$ defined by
$$\sigma_{x_{0}}(s)=\exp_{x_{0}}(sv), \, \, \sigma_{y_{0}}(s)=\exp_{y_{0}}(sw).$$
We want to calculate
$$
d^{2}\varphi(x_{0},y_{0})(v,w)^{2}=
\frac{d^{2}}{ds^{2}}\varphi(\sigma_{x_{0}}(s),
\sigma_{y_{0}}(s))|_{s=0},
$$
where $\varphi(x,y)=d(x,y)^{2}$. To this end let us denote by
$\alpha_{s}:[0,\ell]\to M$ the unique minimizing geodesic joining
the point $\sigma_{x_{0}}(s)$ to the point $\sigma_{y_{0}}(s)$
(now, for $s\neq 0$, $\alpha_{s}$ is not necessarily parameterized
by arc-length), and let us define $\alpha:[0,\ell]\times
[-\varepsilon, \varepsilon]\to M$ by $\alpha(t,s)=\alpha_{s}(t)$.
Then $\alpha$ is a smooth variation through geodesics of
$\gamma(t)=\alpha(t,0)$ and, according to the above discussion, we
have
    $$
    \varphi(\sigma_{x_{0}}(s), \sigma_{y_{0}}(s))=L(s)^{2}=\ell
    E(s),
    $$
and therefore
    $$
    d^{2}\varphi(x_{0},y_{0})(v,w)^{2}=\ell E''(0). \eqno(4)
    $$
If we denote $X(t)=\partial\alpha(t,0)/\partial s$, the
variational field of $\alpha$, then the formula for the second
variation of energy (see \cite[p. 197]{dC}) tells us that
    $$
    \frac{1}{2}E''(0)=-\int_{0}^{\ell}\langle X,
    X''+R(\gamma',X)\gamma'\rangle dt \, + \, \langle X(t),
    X'(t)\rangle |_{t=0}^{t=\ell} \, + \,
    \langle\frac{D}{ds}\frac{\partial\alpha}{\partial s}(t,0),
    \gamma'(t)\rangle |_{t=0}^{t=\ell}, \eqno(5)
    $$
or equivalently
    $$
    \frac{1}{2}E''(0)=\int_{0}^{\ell}\left( \langle X', X'\rangle -
    \langle R(\gamma',X)\gamma', X\rangle \right)dt \, + \,
    \langle\frac{D}{ds}\frac{\partial\alpha}{\partial s}(t,0),
    \gamma'(t)\rangle |_{t=0}^{t=\ell}, \eqno(6)
    $$
where we denote $X'=\nabla_{\gamma'(t)}X$, and
$X''=\nabla_{\gamma'(t)} X'$.

Note that, since the variation field of a variation through
geodesics is always a Jacobi field, and since the points $x_{0}$
and $y_{0}$ are not conjugate, the field $X$ is in fact the unique
Jacobi field along $\gamma$ satisfying that $X(0)=v$, $X(\ell)=w$,
that is $X$ is the unique vector field along $\gamma$ satisfying
    $$
    X''(t)+R(\gamma'(t), X(t))\gamma'(t)=0, \,\,\, \textrm{ and } \,\,\, X(0)=v, \,\,\
    X(\ell)=w,
    $$
where $R$ is the curvature of $M$. On the other hand, since the
curves $s\to\alpha(0,s)=\sigma_{x_{0}}(s)$ and $s\to\alpha(\ell,
s)=\sigma_{y_{0}}(s)$ are geodesics, we have that
    $$
    \langle\frac{D}{ds}\frac{\partial\alpha}{\partial s}(t,0),
    \gamma'(t)\rangle |_{t=0}^{t=\ell}=0.
    $$
These observations allow us to simplify $(5)$ and $(6)$ by
dropping the terms that vanish, thus obtaining that
    $$
    \frac{1}{2}E''(0)=\langle X(\ell),
    X'(\ell)\rangle - \langle X(0),
    X'(0)\rangle \eqno(7)
    $$
and also
    $$
    \frac{1}{2}E''(0)=\int_{0}^{\ell}\left( \langle X', X'\rangle -
    \langle R(\gamma',X)\gamma', X\rangle \right)dt. \eqno(8)
    $$
Recall that the right-hand side of $(8)$ is called the index form
and is denoted by $I(X,X)$.

By combining $(4), (7)$ and $(8)$ we get
\begin{eqnarray*}
& & d^{2}\varphi(x_{0},y_{0})(v,w)^{2}=2\ell \left( \langle
X(\ell),
    X'(\ell)\rangle - \langle X(0),
    X'(0)\rangle \right)=\\
& &= 2\ell \int_{0}^{\ell}\left( \langle X', X'\rangle -
    \langle R(\gamma',X)\gamma', X\rangle \right)dt.
\end{eqnarray*}

Therefore condition $(\sharp)$ holds if and only if, for every
Jacobi field $X$ along $\gamma$ with $X(0)=v,
X(\ell)=w=L_{x_{0}y_{0}}v$, one has that
$$\langle X(\ell),
    X'(\ell)\rangle - \langle X(0),
    X'(0)\rangle\leq 0 \eqno(\diamondsuit)$$
or, equivalently,
    $$\int_{0}^{\ell}\left( \langle X', X'\rangle -
    \langle R(\gamma',X)\gamma', X\rangle \right)dt\leq 0
    $$
for the same Jacobi fields.

\medskip

\noindent {\bf Proof of Proposition \ref{key property of hessian
of distance}:} The proof of $(2)$ is immediate and is well
referenced (see for instance \cite[Theorem IX.4.3]{Lang}): if $M$
has nonpositive sectional curvature then we have $\langle
R(\gamma',X)\gamma', X\rangle\leq 0$, hence, according to the
above formulas,
    $$
    d^{2}\varphi(x_{0},y_{0})(v,w)^{2}= 2\ell \int_{0}^{\ell}\left( \langle X', X'\rangle -
    \langle R(\gamma',X)\gamma', X\rangle \right)dt\geq
    2\ell \int_{0}^{\ell}\langle X', X'\rangle dt\geq 0,
    $$
which proves $(2)$. Note that in this case we do not use that
$w=L_{x_{0}y_{0}}v$, so this holds for all $v, w$.

\medskip

Our proof of $(1)$ uses the following Lemma, which is a
restatement of Corollary 10 in Chapter 8 of \cite{Spivak}.

%%%%%%%%%%%%%%%%%%%
\begin{lem}\label{Spivaks lemma}
Let $\gamma:[0,\ell]\to M$ be a geodesic without conjugate points,
$X$ a Jacobi field along $\gamma$, and $Z$ a piecewise smooth
vector field along $\gamma$ such that $X(0)=Z(0)$ and
$X(\ell)=Z(\ell)$. Then
    $$
    I(X, X)\leq I(Z,Z),
    $$
and equality holds only when $Z=X$.

\noindent {\em That is, among all vector fields along $\gamma$
with the same boundary conditions, the unique Jacobi field along
$\gamma$ determined by those conditions minimizes the index form.
Recall that
$$I(Z,Z)=\int_{0}^{\ell}\left( \langle Z', Z'\rangle -
    \langle R(\gamma',Z)\gamma', Z\rangle \right)dt,$$
but this number is not equal to $\langle Z(\ell),
    Z'(\ell)\rangle - \langle Z(0),
    Z'(0)\rangle$ unless $Z$ is a Jacobi field.}
\end{lem}
%%%%%%%%%%%%%%%%%%%

Let $X$ be the unique Jacobi field with $X(0)=v,
X(\ell)=w=L_{x_{0}y_{0}}(v)$. Define $Z=P(t)$, where $P(t)$ is the
parallel translation along $\gamma$ with $P(0)=v$ (hence
$P(\ell)=w$). The field $Z$ is not necessarily a Jacobi field, but
it has the considerable advantage that $Z'(t)=0$ for all $t$, so
we have that
    $$
    I(Z,Z)=\int_{0}^{\ell}\left( \langle Z', Z'\rangle -
    \langle R(\gamma',Z)\gamma', Z\rangle \right) dt=
    -\int_{0}^{\ell}
    \langle R(\gamma',Z)\gamma', Z\rangle dt\leq 0
    $$
because $M$ has nonnegative sectional curvature. We then deduce
from the above Lemma that
    $$
    \langle X(\ell),
    X'(\ell)\rangle - \langle X(0),
    X'(0)\rangle = I(X,X)\leq I(Z,Z)\leq 0,
    $$
which, according to the above remarks (see $(\diamondsuit)$),
concludes the proof. \hspace{1.5cm} $\Box$
%%%%%%%%%%%%%%%%

\bigskip

Even though we will not have $d^{2}\varphi(x,y)(v,
L_{xy}v)^{2}\leq 0$ when $M$ has negative curvature, we can
estimate this quantity and show that it is bounded by a term of
the order of $d(x,y)^{2}$, provided that the curvature is bounded
below. This will also be used in the next section to deduce a
comparison result which holds for all Riemannian manifolds
(assuming that $F$ is uniformly continuous).

%%%%%%%%%%%%%%%%%%%%
\begin{prop}\label{bound for A with no restriction on curvature}
Let $M$ be a Riemannian manifold. Consider the function
$\varphi(x,y) = d(x,y)^{2}$, defined on $M\times M$. Assume that
the sectional curvature $K$ of $M$ is bounded below, say $K\geq
-K_{0}$. Then
    $$
    d^{2}\varphi(x,y)(v, L_{xy}v)^{2}\leq 2 K_{0}d(x,y)^{2} \|v\|^{2}
    $$
for all $v\in TM_{x}$ and $x,y\in M$ with $d(x,y)<\min\{ i_{M}(x),
i_{M}(y)\}$.
\end{prop}
\noindent Note that for $K_{0}=0$ we recover part $(1)$ of
Proposition \ref{key property of hessian of distance}.
%%%%%%%%%%%%%%%%%%%%
\begin{proof}
Let $X, Z$ be as in the proof of $(1)$ of the preceding
Proposition. With the same notations, we have that
\begin{eqnarray*}
& & d^{2}\varphi(x,y)(v, L_{xy}v)^{2}= 2\ell \, I(X,X)\leq 2\ell
\, I(Z,Z)=\\
& & -2\ell\int_{0}^{\ell}
    \langle R(\gamma',Z)\gamma', Z\rangle dt\leq
    2\ell\int_{0}^{\ell}
    K_{0}|\gamma'(t)\wedge Z(t)|^{2} dt\leq \\
& & 2\ell\int_{0}^{\ell}
    K_{0} \|\gamma'(t)\|^{2} \, \|Z(t)\|^{2}=
    2\ell\int_{0}^{\ell}K_{0}\|v\|^{2}dt=2\ell^{2}K_{0}\|v\|^{2},
\end{eqnarray*}
which proves the result.
\end{proof}
%%%%%%%%%%%%%%%%%%%%

\medskip

%%%%%%%%%%%%%%%%%%%%%%%%%%%%%%%%%%%%%%%%%%%%%%%%%%%%%%%%%%%%%%%%%%%%%%%%%%%%%%%%%

\section{Comparison results for the Dirichlet problem}

In this section and throughout the rest of the paper we will often
abbreviate saying that $u$ is an upper semicontinuous function on
a set $\Omega$ by writing $u\in USC(\Omega)$. Similarly,
$LSC(\Omega)$ will stand for the set of lower semicontinuous
functions on $\Omega$.

The following lemma will be used in the proof of the main
comparison result for the Dirichlet problem
    $$
    F(x,u(x), du(x), d^{2}u(x))=0 \,\, \textrm{ on }
    \,\ \Omega; \,\,\, u=f \,\, \textrm{ on } \,\ \partial\Omega.
    \eqno(DP)
    $$
%%%%%%%%%%%%%%%%%%%%
\begin{lem}\label{alphalemma}
Let $\Omega$ be a subset of a Riemannian manifold $M$, $u\in
USC(\overline{\Omega})$, $v\in LSC(\overline{\Omega})$ and
    $$m_{\alpha}:=\sup_{\Omega\times\Omega}
    (u(x)-v(y)-{\alpha\over 2}d(x,y)^2)$$
for $\alpha>0.$ Let $m_{\alpha}<\infty$ for large $\alpha$ and
$(x_{\alpha},y_{\alpha})$ be such that
    $$\lim_{\alpha\rightarrow \infty}(m_{\alpha}-(u(x_{\alpha})-
    v(y_{\alpha})-{\alpha\over 2}d(x_{\alpha},y_{\alpha})^2))=0.$$
Then we have:
\begin{enumerate}
\item $\lim_{\alpha\rightarrow \infty}\alpha
                  d(x_{\alpha},y_{\alpha})^2=0$, and
\item $\lim_{\alpha\rightarrow\infty}m_{\alpha}=u(\widehat{x})-v(\widehat{x})
=sup_{x\in \Omega}(u(x)-v(x))$ whenever $\widehat{x}\in\Omega$ is
a limit point of $x_{\alpha}$ as $\alpha\rightarrow\infty$.
\end{enumerate}
\end{lem}
%%%%%%%%%%%%%%%%%%%%%%%
\begin{proof}
The result is proved in \cite[Lemma 3.1]{CLI} in the case
$M=\mathbb{R}^{n}$, and the same proof clearly works in the
generality of the statement (in fact this holds in any metric
space).
\end{proof}
%%%%%%%%%%%%%%%%%%%%%%%

Now we can prove the main comparison result for the Dirichlet
problem.

%%%%%%%%%%%%%%%%%%%%%%%%%%%%%%%%
\begin{thm}\label{maximum principle for Dirichlet}
Let $\Omega$ be a bounded open subset of a complete
finite-dimensional Riemannian manifold $M$, and
$F:\mathcal{X}\to\mathbb{R}$ be proper, continuous, and satisfy:
\begin{enumerate}
\item  there exists $\gamma>0$ such that
    $$\gamma(r-s)\leq F(x,r,\zeta,Q)-F(x,s,\zeta,Q)$$
for $r\geq s$; and
\item there exists a function $\omega:[0,\infty]\rightarrow[0,\infty]$
with $\lim_{t\to 0^{+}}\omega(t)=0$ and such that
    $$
    F(y,r, \alpha\exp_{y}^{-1}(x),Q)-F(x,r,-\alpha\exp_{x}^{-1}(y),P)
    \leq \omega\left(\alpha d(x,y)^2+d(x,y)\right)
    $$
for all $x,y\in\Omega$, $r\in\R$, $P\in T_{2, s}(M)_x, Q\in T_{2,
s}(M)_y$ with
\begin{equation*}\label{*}
        -\left(
        \frac{1}{\varepsilon_{\alpha}}+\|A_{\alpha}\|\right)\left(\begin{array}{cc}
 I & 0  \\
 0 & I \\
\end{array}\right)\leq
\left(\begin{array}{cc}
 P & 0  \\
 0 & -Q \\
\end{array}\right)\leq A_{\alpha}+\varepsilon_{\alpha} A_{\alpha}^{2}, \eqno(*)
\end{equation*}
where $A_{\alpha}$ is the second derivative of the function
$\varphi_{\alpha}(x,y)=\frac{\alpha}{2}d(x,y)^{2}$ at the point
$(x,y)\in\ M\times M$,
$$\varepsilon_{\alpha}=\frac{1}{2(1+\|A_{\alpha}\|)}, $$
 and the points $x, y$ are assumed to be
close enough to each other so that $d(x,y)<\min\{ i_{M}(x),
i_{M}(y)\}$.
\end{enumerate}
Let $u\in USC(\overline{\Omega})$ be a subsolution and $v\in
LSC(\overline{\Omega})$ a supersolution of $F=0$ on $\Omega$, and
$u\leq v$ on $\partial \Omega$.

Then $u\leq v$ holds on all of $\overline{\Omega}$.

In particular the Dirichlet problem $(DP)$ has at most one
viscosity solution.
\end{thm}
%%%%%%%%%%%%%%%%%%%%%%%%%%%%%%%%
\begin{proof}
Assume to the contrary that there exists $z\in\Omega$ with
$u(z)>v(z)$.  By compactness of $\overline{\Omega}$ and upper
semicontinuity of $u-v$, and according to Lemma \ref{alphalemma},
there exist $x_\alpha, y_\alpha$ so that, with the the notation of
Lemma \ref{alphalemma},
    $$u(x_{\alpha})-v(y_{\alpha})-{\alpha\over
2}d(x_{\alpha}, y_{\alpha})^2=m_{\alpha}\geq \delta:=u(z)-v(z)>0,
\eqno(3)
$$ and
    $$
    \alpha d(x_{\alpha}, y_{\alpha})^{2}\to 0 \textrm{ as }
    \alpha\to\infty. \eqno(4)
    $$
Again by compactness of $\overline{\Omega}$ we can assume that a
subsequence of $(x_{\alpha}, y_{\alpha})$, which we will still
denote $(x_{\alpha}, y_{\alpha})$ (and suppose $\alpha\in\N$),
converges to a point $(x_0, y_0)\in \overline{\Omega}\times
\overline{\Omega}$. By Lemma \ref{alphalemma} we have that
$x_0=y_0$ and
    $$
    \delta\leq \lim_{\alpha\to\infty}m_{\alpha}=u(x_{0})-v(x_{0})=\sup_{\overline{\Omega}}(u(x)-v(x)),
    $$
and in view of the condition $u\leq v$ on $\partial\Omega$ we have
that $x_{0}\in\Omega$, and $x_{\alpha}, y_{\alpha}\in\Omega$ for
large $\alpha$.

Fix $r_{0}>0$ and $R_{0}>0$ such that, for every $x\in B(x_{0},
r_{0})$, $\exp_{x}$ is a diffeomorphism from $B(0,R_{0})\subset
TM_x$ onto $B(x,R_{0})\supset B(x_{0}, r_{0})$ (see \cite[Theorem
3.7 of Chapter 3]{dC}). Then, for every $x, y\in B(x_{0}, r_{0})$
we have that $d(x,y)<\min\{i_{M}(x), i_{M}(y)\}$,  the vectors
$\exp_{x}^{-1}(y)\in TM_{x}\equiv TM^{*}_{x}$ and
$\exp_{y}^{-1}(x)\in TM_{y}\equiv TM^{*}_{y}$ are well defined,
and the function $\varphi(x,y)=d(x,y)^{2}$ is $C^{2}$ smooth on
$B(x_{0}, r_{0})\times B(x_{0}, r_{0})\in M\times M$. Taking a
subsequence if necessary, we can assume that $x_{\alpha},
y_{\alpha}\in B(x_{0}, r_{0})$ for all $\alpha$.

\medskip

Now, for each $\alpha$, we can apply Theorem \ref{key to
comparison} with $\Omega_{1}=\Omega_{2}=B(x_{0}, r_{0})$,
$u_{1}=u$, $u_{2}=-v$,
$\varphi(x,y)=\varphi_{\alpha}(x,y):=\frac{\alpha}{2}d(x,y)^{2}$,
and for
$$\varepsilon=\varepsilon_{\alpha}:=\frac{1}{2\left(1+\|d^{2}\varphi_{\alpha}(x_{\alpha},
y_{\alpha})\|\right)}.$$ Since $(x_{\alpha},y_{\alpha})$ is a
local maximum of the function $(x,y)\mapsto
u(x)-v(y)-\varphi(x,y)$, we obtain bilinear forms $P\in
\mathcal{L}^{2}_{s}((TM_{x_{\alpha}}, \mathbb{R})$, and $Q\in
\mathcal{L}^{2}_{s}((TM_{y_{\alpha}}, \mathbb{R})$ such that
$$
\left(\frac{\partial}{\partial x}\varphi(x_{\alpha},y_{\alpha}), P \right)\in \overline{J}^{\, 2, +}u(x_{\alpha}),
$$
$$
\left(-\frac{\partial}{\partial y}\varphi(x_{\alpha},y_{\alpha}),
Q \right)\in \overline{J}^{\, 2, -}v(y_{\alpha})
$$
(recall that
$\overline{J}^{2,-}v(y_{\alpha})=-\overline{J}^{2,+}(-v)(y_{\alpha})$),
and
    $$-\left({1\over\varepsilon_{\alpha}}+\|A_{\alpha}\|\right)I\leq
    \left(
\begin{array}{ccc}
  P   &   0    \\
   0  &  -Q  \\
\end{array}\right)\leq A_{\alpha}+\varepsilon_{\alpha} A_{\alpha}^2,$$
where $A_{\alpha}=d^2\varphi(x_{\alpha}, y_{\alpha})\in
\mathcal{L}^{2}_{s}(TM_{(x,y)},\mathbb{R})$, so we get that
condition $(*)$ holds for $x=x_{\alpha}, y=y_{\alpha}$. Therefore,
according to condition $(2)$, we have that
    $$
    F(y_{\alpha},r, \alpha\exp_{y_{\alpha}}^{-1}(x_{\alpha}),Q)-F(x_{\alpha},r,-\alpha\exp_{x_{\alpha}}^{-1}(y_{\alpha}),P)
    \leq \omega(\alpha
    d(x_{\alpha},y_{\alpha})^2+d(x_{\alpha},y_{\alpha})) \eqno(5)
    $$
On the other hand, from equation $(1)$ in the preceding section we
have that
    $$
    \frac{\partial}{\partial x}\varphi(x_{\alpha},y_{\alpha})=-\alpha\exp_{x_{\alpha}}^{-1}(y_{\alpha}),
    \, \textrm{ and } \,
    -\frac{\partial}{\partial
    y}\varphi(x_{\alpha},y_{\alpha})=\alpha\exp_{y_{\alpha}}^{-1}(x_{\alpha}),
    $$
hence $(-\alpha\exp_{x_{\alpha}}^{-1}(y_{\alpha}), P)\in
\overline{J}^{2,+}u(x_{\alpha})$,
$(\alpha\exp_{y_{\alpha}}^{-1}(x_{\alpha}), Q)\in
\overline{J}^{2,-}v(y_{\alpha})$. Since $u$ is subsolution and $v$
is supersolution, and $F$ is continuous we then have, according to
Remark \ref{continuity of F extends inequalities to closures of
subjets}, that
    $$
    F(x_{\alpha}, u(x_{\alpha}), -\alpha\exp_{x_{\alpha}}^{-1}(y_{\alpha}), P)\leq 0\leq
    F(y_{\alpha}, v(y_{\alpha}), \alpha\exp_{y_{\alpha}}^{-1}(x_{\alpha}),
    Q). \eqno(6)
    $$
By combining equations $(3)$, $(4)$, $(5)$ and $(6)$ above, and
using condition $(1)$ too, we finally get
\begin{eqnarray*}
& & 0<\gamma\delta\leq
\gamma\left(u(x_{\alpha})-v(y_{\alpha})\right)\leq\\
& & F(x_{\alpha}, u(x_{\alpha}),
-\alpha\exp_{x_{\alpha}}^{-1}(y_{\alpha}), P)- F(x_{\alpha},
v(y_{\alpha}), -\alpha\exp_{x_{\alpha}}^{-1}(y_{\alpha}), P)\leq \\
& &F(x_{\alpha}, u(x_{\alpha}),
-\alpha\exp_{x_{\alpha}}^{-1}(y_{\alpha}), P) -
    F(y_{\alpha}, v(y_{\alpha}), \alpha\exp_{y_{\alpha}}^{-1}(x_{\alpha}),
    Q) +\\
& & +F(y_{\alpha}, v(y_{\alpha}),
\alpha\exp_{y_{\alpha}}^{-1}(x_{\alpha}),Q)-
F(x_{\alpha},v(y_{\alpha}),-\alpha\exp_{x_{\alpha}}^{-1}(y_{\alpha}),P)\leq
\\ & &\leq \omega(\alpha d(x_{\alpha},y_{\alpha})^{2}+d(x_{\alpha},y_{\alpha})),
\end{eqnarray*}
and the contradiction follows by letting $\alpha\to\infty$.
\end{proof}
%%%%%%%%%%%%%%%%%%%%%%%%%
\begin{rem}
{\em Observe that, since
$\alpha\exp_{y}^{-1}(x)=L_{xy}(-\alpha\exp_{x}^{-1}(y))$,
condition $(2)$ of Theorem \ref{maximum principle for Dirichlet}
can be replaced with a stronger but simpler assumption, namely
that
    $$
    F(y,r, L_{xy}\zeta, Q)-F(x,r,\zeta,P)
    \leq \omega\left(\alpha d(x,y)^2+d(x,y)\right)
    $$
for all $x,y\in\Omega$, $r\in\R$, $P\in T_{2, s}(M)_x, Q\in T_{2,
s}(M)_y$, $\zeta\in TM^{*}_{x}$ satisfying $(*)$.}
\end{rem}
%%%%%%%%%%%%%%%%%%%%%%%%%
\begin{rem}\label{when u and v are bounded r and s can be taken to be bounded}
{\em If we want to compare two solutions $u$ and $v$ of $F=0$ and
we know that these functions are bounded by some $R>0$ (e.g. when
$M$ is compact) then it is obvious from the above proof that it
suffices to require that conditions $(1)$ and $(2)$ of Theorem
\ref{maximum principle for Dirichlet} be satisfied for all $r, s$
in the interval $[-R, R]$.}
\end{rem}
%%%%%%%%%%%%%%%%%%%%%%%%%
\begin{prop}\label{observations on maximum principle for positive curvature}
If $M$ has nonnegative sectional curvature, then condition $(*)$
implies that $P\leq L_{yx}(Q)$.
\end{prop}
%%%%%%%%%%%%%%%%%%%%%%%%%
\begin{proof}
Let $\lambda_{1}, ..., \lambda_{n}$ be the eigenvalues of the
restriction of $A_{\alpha}$ to the subspace $\mathcal{D}=\{(v,
L_{xy}v) : v\in TM_{x}\}$ of $TM_{x}\times TM_{y}$. By Proposition
\ref{key property of hessian of distance} we have that
$A_{\alpha}(v,L_{xy}v)^{2}\leq 0$ for all $v\in TM_{x}$, that is
${(A_{\alpha})}_{|\mathcal{D}}\leq 0$, or equivalently
$\lambda_{i}\leq 0$ for $i=1, ..., n$. With our choice of
$\varepsilon_{\alpha}$, this implies that
    $$
    \lambda_{i}+\varepsilon_{\alpha}\lambda_{i}^{2}\leq
    \lambda_{i}+\frac{1}{2(1+\sup_{1\leq j\leq n}|\lambda_{j}|)}\lambda_{i}^{2}\leq
    \lambda_{i}+\frac{|\lambda_{i}|}{2}=\frac{\lambda_{i}}{2}\leq
    0
    $$
and since $\lambda_{i}+\varepsilon_{\alpha}\lambda_{i}^{2}$, $i=1,
..., n$, are the eigenvalues of
$\left(A_{\alpha}+\varepsilon_{\alpha}A_{\alpha}^{2}\right)_{|\mathcal{D}}$,
this means that
$$
\left(A_{\alpha}+\varepsilon_{\alpha}A_{\alpha}^{2}\right)(v,L_{xy}v)^{2}\leq
0.
$$
Then condition $(*)$ implies that
    $$
    P(v)^{2}-Q(L_{xy}v)^{2}\leq
    \left(A_{\alpha}+\varepsilon_{\alpha}A_{\alpha}^{2}\right)(v,L_{xy}v)^{2}\leq 0
    $$
for all $v\in TM_{x}$, which means that $P\leq L_{yx}(Q)$.
\end{proof}
%%%%%%%%%%%%%%%%%%%

\medskip

Therefore, {\em if $M$ has nonnegative curvature and $F$ is
degenerate elliptic} then $(*)$ automatically implies that
    $$
    F(x,r,\zeta, L_{yx}Q) - F(x,r,\zeta, P)\leq 0,
    $$
hence
\begin{eqnarray*}
& & F(y,r, L_{xy}\zeta,Q)-F(x,r,\zeta,P)=\\
& & F(y,r,
L_{xy}\zeta,Q)-F(x,r, \zeta,L_{yx}Q)+F(x,r,\zeta,L_{yx}Q)-F(x,r,\zeta,P)\leq\\
& & F(y,r, L_{xy}\zeta,Q)-F(x,r, \zeta,L_{yx}Q),
\end{eqnarray*}
and we see that condition $(2)$ of the Theorem is satisfied if we
additionally require, for instance, that
    $$
    F(y,r, \eta,Q)-F(x,r, L_{yx}\eta,L_{yx}Q)\leq \omega(d(x,y)).
    \eqno(2\sharp)
    $$

Note that, in the case $M=\mathbb{R}^{n}$ we have
$L_{yx}\eta\equiv\eta$ and $L_{yx}Q\equiv Q$, and condition
$(2\sharp)$ simply means that $F(x, u, du, d^2 u)$ is uniformly
continuous with respect to the variable $x$. Therefore we can
regard condition $(2\sharp)$ as the natural extension to the
Riemannian setting of the Euclidean notion of uniform continuity
of $F$ with respect to $x$. This justifies the following
%%%%%%%%%%%%%
\begin{defn}
{\em We will say that $F:\mathcal{X}\to\mathbb{R}$ is {\em
intrinsically uniformly continuous with respect to the variable $x$}
if condition $(2\sharp)$ above is satisfied.}
\end{defn}
%%%%%%%%%%%%%
\begin{rem}
{\em As we saw in Example \ref{examples of invariant functions}
above, many interesting examples of equations involving
nonincreasing symmetric functions of the eigenvalues of $d^2 u$
(such as the trace and the positive determinant
$\textrm{det}_{+}$) automatically satisfy condition $(2\sharp)$ as
long as they do not depend on $x$. In fact, since the eigenvalues
of $A$ are the same as those of $L_{yx}A L_{xy}$, any function of
the form
    $$
    F(x, r, \zeta, A)=G(r, \|\zeta\|_x, \textrm{ eigenvalues of } A)
    $$
is intrinsically uniformly continuous with respect to $x$.}
\end{rem}

%%%%%%%%%%%%%%%%%%%%
Therefore, for manifolds of nonnegative curvature, we do not need
to impose that $F$ depends on $d^{2}u(x)$ in a uniformly
continuous manner: the assumptions that $F$ is degenerate elliptic
and intrinsically uniformly continuous with respect to $x$ are
sufficient. Let us sum up what we have just shown.
%%%%%%%%%%%%%%%%%%%%
\begin{cor}\label{comparison Dirichlet positive curvature}
Let $\Omega$ be a bounded open subset of a complete
finite-dimensional Riemannian manifold $M$ with nonnegative
sectional curvature, and $F:\mathcal{X}\to\mathbb{R}$ be continuous,
degenerate elliptic, and satisfy:
\begin{enumerate}
\item $F$ is strongly increasing, that is there exists $\gamma>0$ such that, if $r\geq s$ then
    $$\gamma(r-s)\leq F(x,r,\zeta,Q)-F(x,s,\zeta,Q);$$
\item $F$ is intrinsically uniformly continuous with respect to
$x$ (that is there exists a function
$\omega:[0,\infty]\rightarrow[0,\infty]$ with $\lim_{t\to
0^{+}}\omega(t)=0$ and such that $$
    F(y,r,\eta,Q)-F(x,r, L_{yx}\zeta,L_{yx}Q)\leq \omega(d(x,y))
    \eqno(2\sharp)
    $$
for all $x,y, r, \zeta, Q$ with $d(x,y)<\min\{ i_{M}(x),
i_{M}(y)\}$).
\end{enumerate}
Let $u\in USC(\overline{\Omega})$ be a subsolution and $v\in
LSC(\overline{\Omega})$ a supersolution of $F=0$ on $\Omega$, and
$u\leq v$ on $\partial \Omega$.

Then $u\leq v$ on all of $\overline{\Omega}$.

In particular, the Dirichlet problem $(DP)$ has at most one
viscosity solution.
\end{cor}
%%%%%%%%%%%%%%%%%%%%

When $M$ has negative curvature, condition $(*)$ does not imply
$P\leq L_{yx}Q$, and degenerate ellipticity together with
fulfillment of $(2\sharp)$ is not enough to ensure that condition
$(2)$ of Theorem \ref{maximum principle for Dirichlet} is
satisfied. In this case condition $(2)$ of \ref{maximum principle
for Dirichlet} involves kind of a uniform continuity assumption on
the dependence of $F$ with respect to $d^{2}u(x)$. Let us be more
explicit.
%%%%%%%%%%%%%%%%%%%%
\begin{prop}\label{estimation of norm of P-LQ}
Assume that $M$ has sectional curvature bounded below by some
constant $-K_{0}\leq 0$. Then condition $(*)$ in Theorem
\ref{maximum principle for Dirichlet} implies that
$$P-L_{yx}(Q)\leq \frac{3}{2}K_{0} \, \alpha \, d(x,y)^{2}\, I, $$
where $I(v)^{2}=\langle v, v \rangle =\|v\|^{2}$.
\end{prop}
%%%%%%%%%%%%%%%%%%%%
\begin{proof}
We have that $A_{\alpha}=(\alpha/2)d^{2}\varphi(x,y)$, where
$\varphi(x,y)=d(x,y)^{2}$. According to Proposition \ref{bound for
A with no restriction on curvature} we have
    $$
    d^{2}\varphi(x,y)(v, L_{xy}v)^{2}\leq 2 K_{0}d(x,y)^{2} \|v\|^{2}
    $$
for all $v\in TM_{x}$ and $x,y\in M$ with $d(x,y)<\min\{ i_{M}(x),
i_{M}(y)\}$. Therefore
    $$
    A_{\alpha}(v, L_{xy}v)^{2}\leq\alpha K_{0}d(x,y)^{2} \|v\|^{2}.
    $$
This means that the maximum eigenvalue of the restriction of
$A_{\alpha}$ to $\mathcal{D}:=\{(v, L_{xy}v) : v\in TM_{x}\}$,
which we denote $\lambda_{n}$, satisfies
    $$
    \lambda_{n}\leq \alpha K_{0}d(x,y)^{2}.
    $$
If $\lambda_{1}, ..., \lambda_{n}$ are the eigenvalues of
$(A_{\alpha})_{|\mathcal{D}}$ then
$\lambda_{i}+\varepsilon_{\alpha}\lambda_{i}^{2}$, $i=1, ..., n$,
are those of
$\left(A_{\alpha}+\varepsilon_{\alpha}A_{\alpha}^{2}\right)_{|\mathcal{D}}$.
For a given $i=1, ..., n$, if $\lambda_{i}\leq 0$ then
$\lambda_{i}+\varepsilon_{\alpha}\lambda_{i}^{2}\leq 0$ as in the
proof of of Remark \ref{observations on maximum principle for
positive curvature}. In particular, if $\lambda_{n}\leq 0$ then
$\lambda_{j}\leq 0$ for all $j=1, ..., n$, so we get that
$\lambda_{j}+\varepsilon_{\alpha}\lambda_{j}^{2}\leq 0$ for all
$j=1, ..., n$, which means that
    $$
    \left(A_{\alpha}+\varepsilon_{\alpha}A_{\alpha}^{2}\right)(v,
    L_{xy}v)^{2}\leq 0.
    $$
On the other hand, if $\lambda_{n}\geq 0$ then
$\lambda_{n}+\varepsilon_{\alpha}\lambda_{n}^{2}\geq 0$, and
because the function $[0, +\infty)\ni s\mapsto
s+\varepsilon_{\alpha}s^{2}\in [0, +\infty)$ is increasing, the
maximum eigenvalue of
$\left(A_{\alpha}+\varepsilon_{\alpha}A_{\alpha}^{2}\right)_{|\mathcal{D}}$
is precisely $\lambda_{n}+\varepsilon_{\alpha}\lambda_{n}^{2}$.
This means that
    $$
    \left(A_{\alpha}+\varepsilon_{\alpha}A_{\alpha}^{2}\right)(v,
    L_{xy}v)^{2}\leq
    \lambda_{n}+\varepsilon_{\alpha}\lambda_{n}^{2}
    $$
for all $v\in TM_{x}$. Besides we have, by the choice of
$\varepsilon_{\alpha}$, that
    $$
    \lambda_{n}+\varepsilon_{\alpha}\lambda_{n}^{2}\leq
    \lambda_{n}+\frac{1}{2(1+\sup_{1\leq j\leq n}|\lambda_{n}|)}\lambda_{n}^{2}\leq
    \lambda_{n}+\frac{\lambda_{n}}{2}=\frac{3}{2}\lambda_{n},
    $$
hence
    $$
    \left(A_{\alpha}+\varepsilon_{\alpha}A_{\alpha}^{2}\right)(v,
    L_{xy}v)^{2}\leq\frac{3}{2}\lambda_{n}\leq \frac{3}{2}\alpha K_{0}d(x,y)^{2}\|v\|^{2}.
    $$
In any case (no matter what the sign of $\lambda_{n}$ is) we get
that the above inequality holds. Therefore condition $(*)$ implies
    $$
    P(v)^{2}-Q(L_{xy}(v))^{2}\leq \left(A_{\alpha}+\varepsilon_{\alpha}A_{\alpha}^{2}\right)(v,
    L_{xy}v)^{2}\leq \frac{3}{2}K_{0}\, \alpha d(x,y)^{2} \|v\|^{2}.
    $$
\end{proof}
%%%%%%%%%%%%%%%%%%%%

%%%%%%%%%%%%%%%%%%%%
\begin{cor}\label{comparison for uniformly continuous Fs}
Let $M$ be a complete Riemannian manifold (no assumption on
curvature), and $\Omega$ be a bounded open subset of $M$. Suppose
that $F:\mathcal{X}\to\mathbb{R}$ is proper, continuous, and
satisfies the following uniform continuity assumption: for every
$\varepsilon>0$ there exists $\delta>0$ such that
    $$
    d(x,y)\leq\delta, \,\,\  P-L_{yx}(Q)\leq\delta
    I \,
    \implies F(y, r, L_{xy}\zeta,
Q)-F(x, r, \zeta, P)\leq\varepsilon \eqno(2\flat)
    $$
for all $x,y\in M$ with $d(x,y)<i_{\Omega}$, $r\in\mathbb{R}$,
$\zeta\in TM^{*}_{x}$,
$P\in\mathcal{L}^{2}_{s}(TM_{x},\mathbb{R})$, and
$Q\in\mathcal{L}^{2}_{s}(TM_{y},\mathbb{R})$. Assume also that
there is $\gamma>0$ such that
    $$\gamma(r-s)\leq F(x,r,\zeta,Q)-F(x,s,\zeta,Q) \,\, \textrm{ for all } \,\,
    r\geq s.$$
Then there is at most one viscosity solution of the Dirichlet
problem $(DP)$.
\end{cor}
%%%%%%%%%%%%%%%%%%%%%%
\begin{proof}
Since $M$ is complete, we know from the Hopf-Rinow Theorem that
$\overline{\Omega}$ is compact, hence we have that
$i_{\Omega}=\inf_{x\in\overline{\Omega}} i_{M}(x)>0$. Take a
number $r$ with $0<2r< \inf_{x\in\overline{\Omega}} i_{M}(x)$.
Also by compactness of $\overline{\Omega}$, there exists $K_{0}>0$
such that the sectional curvature is bounded below by $-K_{0}$ on
$\overline{\Omega}$. Therefore we have that
$\varphi(x,y)=d(x,y)^{2}$ is $C^\infty$ smooth on the set
$\{(x,y)\in \overline{\Omega}\times \overline{\Omega} :
d(x,y)<r\}$ and, according to the preceding Remark, if $P, Q$
satisfy condition $(*)$ of Theorem \ref{maximum principle for
Dirichlet} we get $P-L_{yx}(Q)\leq \frac{3}{2}K_{0} \, \alpha \,
d(x,y)^{2} \,  I$ whenever $d(x,y)<r$. Then the uniform continuity
assumption on $F$ yields the existence of a function
$\omega:[0,\infty]\rightarrow[0,\infty]$ with $\lim_{t\to
0^{+}}\omega(t)=0$ and such that
    $$
    F(y,r, \alpha\exp_{y}^{-1}(x),Q)-F(x,r,-\alpha\exp_{x}^{-1}(y),P)
    \leq \omega\left(\alpha d(x,y)^2+d(x,y)\right),
    $$
hence the result follows from Theorem \ref{maximum principle for
Dirichlet}.
\end{proof}
%%%%%%%%%%%%%%%%%%%%
\begin{rem}
{\em As is usual with comparison principles, the proof of Theorem
\ref{maximum principle for Dirichlet} can easily be adapted to
show that the viscosity solutions $u$ of the equations $F=0$
depend continuously on $F$. That is, if $u$ is solution of $F=0$
and $v$ is solution of $G=0$, then
    $$
    \sup_{\overline{\Omega}}|u(x)-v(x)|\leq
    \sup_{(x,r,\zeta, A)\in X}|F(x,r,\zeta, A)-G(x,r,\zeta, A)|
    +\sup_{\partial\Omega}|u-v|.
$$
}
\end{rem}
%%%%%%%%%%%%%%%%%%%%

\medskip

%%%%%%%%%%%%%%%%%%%%

\section{Comparison results without boundary conditions}

The same argument as in the proof of Theorem \ref{maximum
principle for Dirichlet}, with some small changes, yields the
following.

%%%%%%%%%%%%%%%%%%%%%%
\begin{thm}\label{main comparison without boundary}
Let $M$ be a connected, complete Riemannian manifold (without
boundary) such that $i(M)>0$, and $F:\mathcal{X}\to\mathbb{R}$ be
proper, continuous, and satisfy assumptions $(1)$ and $(2)$ of
Theorem \ref{maximum principle for Dirichlet}. Let $u$ be a
subsolution, and $v$ a supersolution, of $F=0$. Assume that $u$ and
$v$ are uniformly continuous and $\lim_{x\to\infty}u(x)-v(x)\leq 0$.
Then $u\leq v$ on $M$. In particular, if $M$ is compact, there is at
most one viscosity solution of $F=0$ on $M$.
\end{thm}
%%%%%%%%%%%%%%%%%%%%%%
Uniform continuity and the inequality at infinity guarantee that
the $m_\alpha$ are attained, so the only difference with the proof
of Theorem \ref{maximum principle for Dirichlet} is that now we
cannot assume that $x_{\alpha}$ and $y_{\alpha}$ converge to some
point $x_{0}$, but we do have that $d(x_{\alpha},
y_{\alpha})<i(M)$ for large $\alpha$, hence all the computations
and estimations in the proof of \ref{maximum principle for
Dirichlet} are still valid.
%%%%%%%%%%%%%%%%%%%%%%
\begin{cor}\label{comparison with no boundary and positive curvature}
Let $M$ be a compact Riemannian manifold of nonnegative sectional
curvature, and let $F:\mathcal{X}\to\mathbb{R}$ be continuous,
degenerate elliptic, strongly increasing in $u$, and intrinsically
uniformly continuous with respect to $x$ (that is, $F$ satisfies
conditions $(1-2)$ or Corollary \ref{comparison Dirichlet positive
curvature}. Let $u$ be a subsolution, and $v$ a supersolution of
$F=0$. Then $u\leq v$ on $M$.
\end{cor}
%%%%%%%%%%%%%%%%%%%%%%
\begin{proof}
The same considerations as in Remark \ref{observations on maximum
principle for positive curvature} apply.
\end{proof}
%%%%%%%%%%%%%%%%%%%%%%
\begin{ex}\label{dependence on norm of zeta}
{\em As we remarked above, condition $(2)$ of Corollary
\ref{comparison Dirichlet positive curvature} is easily satisfied
when $F(x,r,\zeta, A)$ does not depend on $\zeta$ and $A$
themselves, but on $\|\zeta\|$ and the eigenvalues of $A$. For
instance, the function
$$F(x,r,\zeta,
A)=r-\left(\textrm{det}_{+}(A)\right)^{3}\|\zeta\|^{2}-f(x)(\textrm{trace}(A))^{5}$$
satisfies $(1)$ and $(2)$ of the above Corollary provided that
$f\geq 0$ and $f$ is uniformly continuous. Therefore the equation
    $$
    u-\left(\textrm{det}_{+}(D^{2}u)\right)^{3}\|\nabla
    u\|^{2}-\left(\Delta u \right)^{5} f =0
    $$
has at most one viscosity solution on any compact manifold of
positive curvature if we only require that $f$ is continuous and
nonnegative.}
\end{ex}
%%%%%%%%%%%%%%%%%%%%%%%%%
\begin{cor}
Let $M$ be a compact Riemannian manifold (no assumption on
curvature). Suppose that $F:\mathcal{X}\to\mathbb{R}$ satisfies the
uniform continuity assumption, and the growth assumption, of
Corollary \ref{comparison for uniformly continuous Fs}. Let $u$ be a
subsolution, and $v$ a supersolution of $F=0$. Then $u\leq v$ on
$M$. In particular there is at most one viscosity solution of $F=0$.
\end{cor}
%%%%%%%%%%%%%%%%%%%%%%
\begin{proof}
Since $M$ is compact the sectional curvature of $M$ is bounded on
all $M$, say $K\geq -K_{0}$. Take a number $r$ with $0<2r<i(M)$.
The function $\varphi(x,y)=d(x,y)^{2}$ is $C^\infty$ on the set
$\{(x,y)\in M\times M : d(x,y)\leq 2r\}$. Suppose that $P$ and $Q$
satisfy $(*)$ of Theorem \ref{maximum principle for Dirichlet}.
Then, from Remark \ref{estimation of norm of P-LQ}, we get that
$P-L_{yx}Q\leq\frac{3}{2} K_{0} \, \alpha \, d(x,y)^{2} \, I$
provided that $d(x,y)<r$. Therefore the uniform continuity
property of $F$ gives us a function
$\omega:[0,\infty]\rightarrow[0,\infty]$ with $\lim_{t\to
0^{+}}\omega(t)=0$ and such that
    $$
    F(y,r, \alpha\exp_{y}^{-1}(x),Q)-F(x,r,-\alpha\exp_{x}^{-1}(y),P)
    \leq \omega\left(\alpha d(x,y)^2+d(x,y)\right).
    $$
Hence we can apply Theorem \ref{main comparison without boundary}
and conclude the result.
\end{proof}
%%%%%%%%%%%%%%%%%%%%%%
%%%%%%%%%%%%%%%%%%%%%%

\medskip

%%%%%%%%%%%%%%%%%%%%%%%%%%%%%%%%%%%%%%%%%%%%%%%%%%%%%%%%%%%%%%%%%%%%%%%%%%%%%%%%%

\section{Existence results}

Perron's method can easily be adapted to the Riemannian setting to
establish existence of viscosity solutions to the Dirichlet
problem. The proof goes exactly as in \cite{CLI} with appropriate
changes. The only step which is not completely obvious is the
proof of the following
%%%%%%%%%%%%%%%%%%%%
\begin{prop}
Let $(\zeta,A)\in J^{2,+}f(z)$ Suppose that $f_n$ is a sequence of
upper semicontinuous functions such that
\begin{enumerate}
\item[{(i)}] there exists $x_n$ such that $(x_n, f_n(x_n))\rightarrow
(x, f(x))$, and
\item[{(ii)}] if  $y_n\rightarrow y$,  then $\limsup_{n\rightarrow\infty}f_n(y_n)\leq
f(y)$.
\end{enumerate}
Then there exist $\widehat{x}_n$ and $(\zeta_n, A_n)\in
J^{2,+}f_n(\widehat{x}_n)$ such that $(\widehat{x}_n
f_n(\widehat{x}_n),\zeta_n, A_n)\rightarrow (x, f(x), \zeta, A)$.
\end{prop}
%%%%%%%%%%%%%%%%%%%%%%
\begin{proof}
Consider the functions $f\circ\exp_{x}$ and $f_{n}\circ\exp_{x}$
defined on a neighborhood of $0$ in $TM_{x}$. These functions
satisfy properties $(i)$ and $(ii)$ of the statement (when they
take the roles of $f$ and $f_{n}$ and $M$ is replaced with
$TM_{x}$). By Corollary \ref{def equivalente por composicion con
exp} we have that $(\zeta, A)\in J^{2,+}(f\circ\exp_{x})(0)$. And
of course the result is known in the case when $M=\mathbb{R}^{n}$,
so we get a sequence $\widehat{v}_{n}$ and
$(\widetilde{\zeta}_{n}, \widetilde{A}_{n})\in
J^{2,+}(f_{n}\circ\exp_{x})(\widehat{v}_{n})$ such that
    $$
    (\widehat{v}_{n}, f_{n}\circ\exp_{x}(\widehat{v}_{n}), \widetilde{\zeta}_{n}, \widetilde{A}_{n})\to
    (0, f\circ\exp_{x}(0), \zeta, A).
    $$
Set $\widehat{x}_{n}=\exp_{x}(\widehat{v}_{n})$. We have that
$\widehat{x}_{n}\to x$ and $f_{n}(x_{n})\to f(x)$. Since
$(\widetilde{\zeta}_{n}, \widetilde{A}_{n})\in
J^{2,+}(f_{n}\circ\exp_{x})(\widehat{v}_{n})$ there exist
functions $\psi_{n}$ such that $f_{n}\circ\exp_{x}-\psi_{n}$
attains a maximum at $\widehat{v}_{n}$,
$\widetilde{\zeta}_{n}=d\psi_{n}(\widehat{v}_{n})$ and
$\widetilde{A}_{n}=d^{2}\psi_{n}(\widehat{v}_{n})$. Let us define
$\varphi=\psi_{n}\circ\exp_{x}^{-1}$ on a neighborhood of $x$.
Then $f_{n}-\varphi_{n}$ attains a maximum at $\widehat{x}_{n}$
so, if we set $\zeta_{n}=d\varphi(\widehat{x}_{n})$,
$A_{n}=d^{2}\varphi(\widehat{x}_{n})$, we have that $(\zeta_{n},
A_{n})\in J^{2,+}f_{n}(\widehat{x}_{n})$. It only remains to show
that $\zeta_{n}\to\zeta$ and $A_{n}\to A$. But this is exactly
what was shown in $(\Leftarrow)$ of the proof of Proposition
\ref{closures of subjets}.
\end{proof}
%%%%%%%%%%%%%%%%%%%%%%

\medskip

By using this Proposition one can prove, as in \cite{CLI},
existence of viscosity solutions to the Dirichlet problem
    $$
    F(x,u, du, d^{2}u)=0 \,\, \textrm{ in } \Omega, \,\,
    \,\,\, u=f \textrm{ on } \partial\Omega, \eqno(DP)
    $$
where $\Omega$ is an open bounded subset of a complete Riemannian
manifold $M$.
%%%%%%%%%%%%%%%%%%%%%%%
\begin{thm}\label{Perron}
Let comparison hold for $(DP)$, i.e., if $w$ is a subsolution of
$(DP)$ and $v$ is a supersolution of $(DP)$, then $w\leq v$.
Suppose also that there exists a subsolution $\underline{u}$ and a
supersolution $\overline{u}$ of $(DP)$ that satisfy the boundary
condition $\underline{u}_*(x)=\overline{u}^*(x)=f(x)$ for $x\in
\partial\Omega.$ Then
\begin{equation*}
  W(x)=sup\{w(x):\underline{u}\leq w\leq \overline{u} \text{ and }
    w \text{ is a subsolution of $(DP)$}\}
\end{equation*}
is a solution of $(DP)$.
\end{thm}
%%%%%%%%%%%%%%%%%%%%%%
\noindent Here we used the following notation:
\begin{equation*}
\begin{array}{c}
 u^*(x)=\lim_{r\downarrow 0}\sup\{u(y):y\in \Omega\text{ and }d(y,x)\leq r\}; \\
 u_*(x)=\lim_{r\downarrow 0}\inf\{u(y):y\in\Omega\text{ and }d(y,x)\leq r\},  \\
\end{array}
\end{equation*}
that is $u^*$ denotes the upper semicontinuous envelope of $u$
(the smallest upper semicontinuous function, with values in
$[-\infty,\infty]$, satisfying $u\leq u^*$), and similarly $u_*$
stands for the lower semicontinuous envelope of $u$.

\medskip

One can also easily adapt the proof of \cite[Theorem 6.17]{AFL2}
to the second order situation, obtaining the following.
%%%%%%%%%%%%%%%%%%%%%%
\begin{cor}
Let $M$ be a compact Riemannian manifold, and $G(x, du, d^2u)$ be
degenerate elliptic and uniformly continuous in the sense of
Corollary \ref{comparison for uniformly continuous Fs}. Then there
exists a unique viscosity solution of the equation \,\,\, $u+G(x,
du, d^{2}u)=0$ on $M$.
\end{cor}
%%%%%%%%%%%%%%%%%%%%%%
\noindent Again, if $M$ has nonnegative curvature, the assumptions
that $F$ is elliptic and intrinsically uniformly continuous with
respect to $x$ are sufficient in order to get an analogous result.

%%%%%%%%%%%%%%%%%%%%%%%%

\medskip

\section{Examples}

Most of the examples of proper $F$'s given in \cite{CLI} remain
valid in the Riemannian setting. In particular, as we have already
seen, the functions $(x,r,\zeta,A)\mapsto -\textrm{det}_{+}(A)$
and $(x,r,\zeta,A)\mapsto -\textrm{trace}(A)$ are degenerate
elliptic and intrinsically uniformly continuous with respect to
$x$. The same is true of all many symmetric functions of the
eigenvalues of $A$, such as minus the minimum (or the maximum)
eigenvalue, and of course nondecreasing combinations and sums of
these are degenerate elliptic too. One can find lots of examples
of nonlinear equations for which the results of this paper yield
existence and uniqueness of viscosity solutions. For instance, one
can easily show that, for every compact manifold of positive
curvature, the equation
    $$
    \max\{u-\lambda_{1}(D^{2}u)\|\nabla u\|^{p}-
    (\Delta u)^{2q+1}\|\nabla u\|^{r}-\left(\textrm{det}_{+}(D^{2}u)\right)^{2k+1} f^{2}, \, \, u-g
    \}=0
    $$
(where $\lambda_{1}$ denotes the minimum eigenvalue function and
$p, q, r, k\in\mathbb{N}$) has a unique viscosity solution if we
only require that $f$ and $g$ are continuous. This gives an idea
of the generality of the above results.

Of course this example is rather unnatural. Let us finish this
paper by examining what our results yield in the case of a classic
equation, that of Yamabe's, which has been extensively studied and
completely solved by using variational methods. We do not claim
that the following discussion gives any new insight into Yamabe's
problem, we only want to study, from the point of view of the
viscosity solutions theory, a well known example of a nonlinear
equation arising from an important geometrical problem.

%%%%%%%%%%%%%%%%%%%%%%
\begin{ex}[The Yamabe equation]\label{Yamabe}
{\em A fundamental problem in conformal geometry is to know
whether or not there exists a conformal metric $g'$ with constant
scalar curvature $S'$ on a given compact $n$-dimensional
Riemannian manifold $(M,g)$, with $n\geq 3$, see \cite{Aubin,
Viaclovsky}. This is equivalent to solving the equation
    $$
    -4\frac{n-1}{n-2} \Delta u +S(x) u= S'
    u^{\frac{n+2}{n-2}},\eqno(Y)
    $$
where $S$ is the scalar curvature of $g$. One can write this
equation in the form $F=0$, where
    $$F(x,r,\zeta, A)=S(x) r -S' r^{\frac{n+2}{n-2}} - 4\frac{n-1}{n-2}\textrm{trace}(A) =0.
    $$
It is clear that $F$ is degenerate elliptic. Assume that $S$ is
everywhere positive and that $S'\leq 0$. Then, by compactness,
there exists $\gamma>0$ such that $S(x)\geq \gamma$ for all $x\in
M$. According to Remark \ref{when u and v are bounded r and s can
be taken to be bounded}, in order to check conditions $(1)$ and
$(2)$ of Theorem \ref{main comparison without boundary} we may
assume that $r, s$ lie on a bounded interval. We have that
    $$
    F(y,r, \eta, Q)-F(x,r, L_{yx}\eta, L_{yx}Q)\leq r| S(y)-S(x)|,
    $$
hence, because $S$ is uniformly continuous on $M$ and $r$ is
bounded, we deduce that $F$ satisfies $(2)$ of Corollary
\ref{comparison with no boundary and positive curvature}. On the
other hand, if $r\geq s$ then
    $$
    F(x, r, \zeta, A)-F(x, s, \zeta, A)=S(x)(r-s)-S'(r^{\frac{n+2}{n-2}}-s^{\frac{n+2}{n-2}})\geq
    \gamma (r-s),
    $$
so condition $(1)$ is also satisfied. It follows that there is at
most one viscosity solution of $F=0$. Existence can be shown by
using Perron's method. In all, we see that if $S$ is everywhere
positive and $S'\leq 0$ then there exists a unique viscosity
solution $u$ of $(Y)$. }
\end{ex}
%%%%%%%%%%%%%%%%%%%%%%

\medskip

%%%%%%%%%%%%%%%%%%%%%%%%%%%%%%%%%
%%%%%%%%%%%%%%%%%%%%%%%%%%%%%%%%%
%%%%%%%%%%%%%%%%%%%%%%%%%%%%%%%%%

\end{document}